\documentstyle[11pt, amssymb, amscd]{article}
\newtheorem{thm}{Theorem}[section]
\newtheorem{lem}[thm]{Lemma}
\newtheorem{prop}[thm]{Proposition}

\newtheorem{cor}[thm]{Corollary}
\newtheorem{defi}[thm]{Definition}

\newcommand{\inv}{^{-1}}
\newcommand{\iso}{\stackrel{\sim}{\longrightarrow}}
\newcommand{\tim}{^\times}
\newcommand{\C}{{\mathbb{C}}}
\newcommand{\Z}{{\mathbb{Z}}}
\newcommand{\R}{{\mathbb{R}}}
\newcommand{\Q}{{\mathbb{Q}}}
\newcommand{\PP}{{\mathbb{P}}}
\newcommand{\Gm}{{\mathbb{G}}_m}

\newcommand{\cH}{{\cal{H}}}

\newcommand{\cI}{{\cal{I}}}
\newcommand{\om}{\Omega}

\newcommand{\cL}{{\cal{L}}}

\newcommand{\ov}{\overline}
\newcommand{\dist}{\mbox{dist}}

\newcommand{\distor}{\mbox{distor}}
\newcommand{\Tr}{\mbox{Tr}}
\newcommand{\ad}{\mbox{ad}}

\newcommand{\Lie}{{\mbox{Lie}\,}}
\newcommand{\diag}{\mbox{diag}}
\newcommand{\ga}{\gamma}
\newcommand{\Ad}{\mbox{Ad}}

\newcommand{\bg}{{\bf  G}}
\newcommand{\bn}{{\bf N}}
\newcommand{\bt}{{\bf T}}
\newcommand{\bu}{{\bf U}}
\newcommand{\bq}{{\bf Q}}
\newcommand{\bp}{{\bf P}}

\begin{document}
\title{Arakelov intersection indices of linear cycles and the geometry of buildings and symmetric spaces}
{ \author{Annette Werner \\ \small Mathematisches Institut, Universit\"at M\"unster, Einsteinstr. 62, D -  48149 M\"unster\\ \small e-mail: werner@math.uni-muenster.de}}
\date{}
\maketitle
\centerline{\bf Abstract}
\small

This paper generalizes Manin's approach towards a geometrical interpretation of Arakelov theory at infinity to linear cycles in projective spaces. 
We show how to interpret  certain non-Archimedean Arakelov intersection numbers of linear cycles on $\PP^{n-1}$ with the combinatorial geometry of the Bruhat-Tits building associated to $PGL(n)$. This geometric setting  has an Archimedean analogue, namely the Riemannian symmetric space associated to $SL(n,\C)$, which we use  to interpret analogous  Archimedean intersection numbers of linear cycles in a similar way.

\normalsize
\vspace*{0.5cm}
\centerline{{\bf MSC} (2000): 14G40, 14M15, 53C35}

\section{Introduction}
In this paper we provide a geometrical interpretation of certain local Arakelov intersection indices of linear cycles on projective spaces. In the non-Archimedean case the corresponding geometric framework is the combinatorial geometry of the Bruhat-Tits building for $PGL$. In the Archimedean case we use the Riemannian geometry of the symmetric space corresponding to $SL(n,\C)$. 

Our motivation was the desire to generalize the results of Manin's paper \cite{ma} to higher dimensions. 

In \cite{ma}, Manin's goal is to enrich the picture of Arakelov theory at the infinite places by constructing a differential-geometric object playing the role of a ``model at infinity''.

He suggests such an object, a certain hyperbolic $3$-manifold, in the case of curves, and he corroborates his suggestion by interpreting various Arakelov intersection numbers in terms of geodesic configurations on this space. 

It is certainly desirable to find such a differential-geometric object also for higher-dimen\-sional varieties, but up to now there have been no results in this direction. One of the goals of this paper is to  present a candidate in the case of projective spaces of arbitrary dimension.

A good strategy for finding such a ``space at infinity'' is to look for a geometric object at the non-Archimedean places which is closely related to a non-Archimedean model und which has an Archimedean analogue.

The idea here is to consider the Bruhat-Tits building $X$ for the group $G = PGL(V)$, where $V$ is an $n$-dimensional vector space over a non-Archimedean local field $K$ of characteristic $0$. The vertices in $X$ correspond to the homothety classes $\{M\}$ of $R$-lattices $M$ in $V$, where $R$ is the ring of integers in $V$.

The boundary $X_\infty$ of $X$ in the Borel-Serre compactification can be identified with the Tits building of $G$, which is just the flag complex in $V$. Thereby the vertices in $X_\infty$ correspond to the non-trivial subspaces of $V$. We write $\PP(W)$ for the irreducible subscheme of the projective space $\PP(V)$ induced by a linear subspace $W \subset V$, and we show that for any fixed vertex $v = \{M\}$ in $X$ the half-geodesic $[v, W]$ connecting  $v$ with the boundary point induced by $W$ governs the reduction of $\PP(W)$ in the model $\PP(M)$ of $\PP(V)$ in the following way: $[v,W]$ and $[v,W']$ share the first $m +1$ vertices iff the reductions of the closures of $\PP(W)$ and $\PP(W')$ in $\PP(M)$ modulo $\pi^m$ (where $\pi$ is a prime element in $R$) coincide.

Then we show how to express a certain intersection index of homologically trivial linear cycles with the combinatorial geometry of $X$. Consider subspaces $A$ and $B$ of dimension $p$ and $C$ and  $D$ of dimension $q = n-p$ of $V$, such that the cycles $\PP(A) - \PP(B)$ and $\PP(C) - \PP(D)$ on $\PP(V)$ have disjoint supports. Then the local Arakelov intersection number of the closures of these cycles in $\PP(M)$ is defined and independent of the choice of a lattice $M$ in $V$. We denote it by $<\PP(A)-\PP(B) , \PP(C)-\PP(D)>$. Under certain conditions (which are e.g. fulfilled if $p=1$ and the intersection is non-trivial), we define explicitely an oriented  geodesic $\gamma$ such that
\[ <\PP(A) - \PP(B), \PP(C) - \PP(D)>= p\;  \distor_\gamma (A \ast \gamma, B \ast \gamma),\]
where $A \ast \gamma$ is the point on $\gamma $ closest to the boundary point induced by $A$ in a suitable sense, and where $\distor$ means oriented distance along the oriented geodesic $\gamma$. 

In \cite{we}, we show another result in this direction. Namely, we give a geometrical interpretation of the intersection index of several arbitrary linear cycles meeting properly on some model $\PP(M)$. So in fact, the geometrical interpretation of non-Archimedean intersections can be pushed quite far. 

The building $X$ has an Archimedean analogue, namely the symmetric space $Z$ corresponding to $SL(n,\C)$. We also have a compactification of $Z$ sharing many features with the Borel-Serre compactification of $X$.
It can be obtained by attaching to $Z$ the set $Z(\infty)$ of geodesic rays in $Z$ emanating at a fixed point $z$ (see e.g. \cite{bgs}). Similar to the Borel-Serre boundary in the non-Archimedean case, $Z(\infty)$ has a decomposition into faces, so that the partially ordered set of faces can be identified with the partially ordered set of proper parabolic subgroups of $SL(n,\C)$. Thereby maximal parabolics correspond to minimal faces, which are points. In this way every non-trivial subspace $W$ of $\C^n$ gives rise to a point in $Z(\infty)$. 

We prove that  geodesics in $Z$ connecting two boundary points corresponding to  subspaces $W$ and $W'$ of $\C^n$ have similar features as geodesics in the building $X$ connecting two boundary points given by subspaces of $V$. More precisely,  there exists a geodesic in $X$ (respectively $Z$) connecting the points corresponding to the subspaces $W$ and $W'$, iff these are complementary in $V$ (respectively $\C^n$). Moreover, the set of geodesics between these points is in bijection with the set of pairs of homothety classes of lattices (respectively hermitian metrics) on $W$ and $W'$. This fits very nicely into Deligne's picture of analogies, where lattices on the non-Archimedean side correspond to hermitian metrics on the Archimedean side (see \cite{de}).

We conclude this paper by interpreting certain Archimedean intersection indices with geodesic configurations in $Z$. Consider $p$-dimensional subspaces $A$ and $B$ and $q = (n-p)$-dimensional subspaces $C$ and $D$ of $\C^n$ such that the cycles $\PP(A) - \PP(B)$ and $\PP(C) - \PP(D)$ on $\PP^{n-1}(\C)$ have disjoint supports. Then we use the Levine currents for these linear cycles  to define their local Archimedean intersection number $<\PP(A) - \PP(B) , \PP(C) -\PP(D)>$. Under certain conditions (which are e.g. fulfilled if $p=1$ and the intersection is nonötrivial), we define explicitely an oriented geodesic $\ga$ such that 
\[ <\PP(A) - \PP(B) , \PP(C) - \PP(D)>= \frac{\sqrt{p}}{\sqrt{q}} \; \distor_\gamma (A \ast \gamma, B \ast \gamma),\]
where $A \ast \gamma$ is again the point on $\gamma $ ``closest'' to the boundary point induced by $A$, namely the orthogonal projection of this point to $\gamma$. (This formula specializes to a formula in \cite{ma}, if $n=2$ and $p = q =1$.)

Hence we get completely parallel formulas in the Archimedean and the non-Archime\-dean picture.

This result and the similar behaviour of geodesics in both cases may suggest to regard $Z$ as some kind of ``model at infinity'' for the projective space
$\PP^{n-1}_{\C}$, and  the set of half-geodesics in $Z$ leading to vertices in $Z(\infty)$  as ``$\infty$-adic reductions'' of linear cycles.

{\bf Acknowledgements: } I would like to thank S. Bloch, L. Br\"ocker, A. Deitmar, Ch. Deninger,  G. Kings, K. K\"unnemann, E. Landvogt, Y. I. Manin, P. Schneider, K. Stramm, M. Strauch and E. de Shalit for useful and inspiring discussions concerning this paper. I am also grateful to the Max-Planck-Institut f\"ur Mathematik in Bonn for financial support and the stimulating atmosphere during the early stages of this work. 

\section{The building and its compactification}
Throughout this paper we denote by $K$ a finite extension of $\Q_p$, by $R$ its valuation ring and by $k$ the residue class field. Besides, $v$ is the valuation map, normalized so that it maps a prime element to $1$. We write $q$ for the cardinality of the residue class field, and we normalize the absolute value on $K$ so that $|x|= q^{-v(x)}$. 

Let $V$ be an $n$-dimensional vector space over $K$. Let us briefly recall the definition of the Bruhat-Tits building $X$ for $\bg = PGL(V)$ (see \cite{brti1} and \cite{la}).

We fix a maximal $K$-split torus $\bt$ and let $\bn = N_G \bt$ be its normalizer. Note that $\bt$ is equal to its centralizer in $\bg$. We write  $G = \bg(K)$, $T=\bt(K)$ and $N= \bn(K)$ for the groups of rational points. 
By $X_\ast(\bt)$ respectively $X^\ast(\bt)$ we denote the cocharacter respectively the character group of $\bt$. 
We have a natural perfect pairing
\begin{eqnarray*}
<,>: & X_\ast(\bt) \times X^\ast(\bt) & \longrightarrow  \Z \\
~ & (\lambda, \chi) & \longmapsto  <\lambda, \chi> ,
\end{eqnarray*}
where $<\lambda, \chi>$ is the integer such that $\chi\circ \lambda(t) = t^{<\lambda, \chi>}$ for all $t \in \Gm$.
Let $\Lambda$ be the $\R$-vector space $\Lambda= X_\ast(\bt) \otimes_\Z \R$. We can identify the dual space $\Lambda^\ast$ with $X^\ast(\bt) \otimes_\Z \R$, and extend $<,>$ to a pairing
\[ <,>: \Lambda \times \Lambda^\ast \longrightarrow \R.\]
Since $<,>$ is perfect, there exists a unique homomorphism $\nu: T \rightarrow \Lambda$ such that
\[ < \nu(z), \chi> = - v(\chi(z))\]
for all $z \in T$ and $\chi \in X^\ast(\bt)$. 
We fix a basis $v_1, \ldots ,v_n$ of $V$ such that 
$ \bt$ is induced by the group of diagonal matrices in $GL(V)$ with respect to  $v_1,\ldots,v_n$. The group $W = N/T$ is the Weyl group of the corresponding root system, hence it acts as a group of reflections on $\Lambda$, and we have a natural homomorphism
$W \longrightarrow GL(\Lambda)$. We can embed $W$ in $N$ as the group of permutation matrices with respect to $v_1,\ldots, v_n$. (A permutation matrix is a matrix which has exactly one entry $1$ in every line and column and which is zero otherwise.) Thereby $N$ is the semidirect product of $T$ and $W$.

Since $\mbox{Aff}(\Lambda) = \Lambda \rtimes GL(\Lambda)$, we can extend $\nu$ to a map
\[ \nu: N = T \rtimes W \longrightarrow \Lambda \rtimes GL(\Lambda) = \mbox{Aff} (\Lambda).\]
The pair $(\Lambda, \nu)$ is called the  empty appartment given by $\bt$ (see \cite{la}, 1.9), and whenever we think of it as an appartment, we  write $A = \Lambda$.

One can define a collection of affine hyperplanes in $A$ decomposing $A$ into infinitely many faces, which are topological simplices (see \cite{la}, \S 11).

Besides,  one defines for every $x \in A$ a certain subgroup $P_x \subset G$ (the would-be stabilizer of x), see \cite{la}, \S 12. Then the building $X$ is given as
\[ X = G \times A /\sim,\]
where the equivalence relation $\sim$ is defined as follows:
\begin{eqnarray*}
(g,x) \sim (h,y) & \mbox{iff  there exists an element } n \in N \\
~ & \mbox{such that } \nu(n)x= y \mbox{ and } g\inv h n \in P_x.\end{eqnarray*}
We have a natural action of $G$ on $X$ via left multiplication on the first factor, and we can embed the appartment $A$ in $X$, mapping $a \in A$ to the class of $(1,a)$. 
This is injective (see \cite{la}, Lemma 13.2).  
For $x \in A$ the group $P_x$ is the stabilizer of $x$. A subset of $X$ of the form $g A$ for some $g \in G$ is called appartment in $X$. Similarly, we define the faces in $g A$ as the subsets $g F$, where $F$ is a face in $A$. Then two points (and even two faces) in $X$ are always contained in a common appartment (\cite{la}, Proposition 13.12 and \cite{brti1}, 7.4.18). Any appartment which contains a point of a face contains the whole face, and even its closure (see \cite{la}, 13.10, 13.11, and \cite{brti1}, 7.4.13, 7.4.14). We fix once and for all a $W$-invariant scalar product on $\Lambda$ which exists by \cite{bou}, VI, 1.1 and 1.2. This induces a metric on $A$. Using the $G$-action it can be continued to a metric $d$ on the whole of $X$ (see \cite{la}, 13.14 and \cite{brti1}, 7.4.20).

We denote by $X^0$ the set of vertices (i.e. $0$-dimensional faces) in $X$. We define a  simplex in $X^0$ to be a subset $\{x_1,\ldots, x_k\}$ of $X^0$ such that $x_1,\ldots,x_k$ are the vertices of a face in $X$. 

Let $\eta_i : \Gm  \rightarrow \bt$ be the cocharacter induced by 
mapping $x$ to the diagonal matrix with diagonal entries $d_1,\ldots,d_n$ such that $d_k= 1$ for $k \neq i$ and $d_i = x$. Then $\eta_1,\ldots, \eta_{n-1}$ is an $\R$-basis of $\Lambda$, and the set of vertices in $A$ is equal to $\bigoplus_{i=1}^{n-1} \Z \eta_i$.

Let $\cL$ be the set of all homothety classes of $R$-lattices of full rank in $V$. We write $\{M\}$ for the class of a lattice $M$. Two different lattice classes $\{M'\}$ and $\{N'\}$ are called adjacent, if there are representatives $M$ and $N$ of $\{M'\}$ and $\{N'\}$ such that 
\[\pi N \subset M \subset N.\]

This relation defines a flag complex, namely  the simplicial complex with vertex set $\cL$ such that the simplices are the sets of pairwise adjacent lattice classes. We have a natural $G$-action on $\cL$ preserving the simplicial structure.

Moreover,  there is a $G$-equivariant bijection 
\[ \varphi: \cL \longrightarrow X^0\]
preserving the simplicial structures (see \cite{we}, section 4). 
If $\{N\} \in \cL$ can be written as $\{N\} = g\{M\}$ for some $g \in G$ and  $M = \pi^{k_1} R v_1 + \ldots + \pi^{k_n} R v_n$, then $\varphi(\{N\})$ is given by the pair $(g, \varphi\{M\}) \in G \times A$, where 
\[\varphi(\{M\}) = \sum_{i=1}^{n-1}  (k_n-k_i) \eta_i\]
is a vertex in $A$.

From now on  we will  identify the vertices in $X$ with $\cL$ without explicitely mentioning the map $\varphi$. 

\begin{defi}The combinatorial distance $\dist(x,y)$ between two points $x$ and $y$ in $X^0$ is defined as 
\begin{eqnarray*}
\dist(x,y) &= &\min \{ k: \mbox{ there are vertices }x= x_0, x_1 ,\ldots, x_k = y,\\
~  &~&\quad \quad \mbox{so that }x_i \mbox{ and }x_{i+1}\mbox{ are adjacent for all }i=0,\ldots, k-1.\}.
\end{eqnarray*}
\end{defi}
Hence $\dist$ is the minimal number of $1$-simplices forming a path between $x$ and $y$. Note that $\dist$ is in general not proportional to the metric $d$ on $X$. 

If $x = \{M\}$ and $y = \{L\}$ are two vertices in $X$, we have $\dist(x,y) = s-r$, where $s= \min\{k: \pi^k L \subset M\}$ and $r =  \max \{k : M \subset \pi^k L \}$ (see \cite{we}, Lemma 4.2). 

Borel and Serre have defined a compactification of  $X$ by attaching the Tits building for $G$ at infinity (see \cite{bose}), which we will now briefly describe.  

First we compactify the appartment $A$. For any half-line $c$ in $A$ and any point $a \in A$ there exists a unique half-line starting in $a$ which is parallel to $c$. We denote it by $[a,c]$. 

Now fix a point $a \in A$, and  $A_\infty$ be the set of halflines in $A$ starting in $a$. Then we define  $\ov{A} = A \cup A_\infty$. 
For all $x \in A$, $c \in A_{\infty}$, and all $\epsilon >0$ we define the cone $C_x(c,\epsilon)$ as
\[C_x(c,\epsilon)=\{z \in \ov{A}: z \neq x \quad \mbox{and} \quad \prec_x([x,c],[x,z]) <\epsilon\}.\]
Here $[x,z]$ is the line from $x$ to $z$ if $z \in A$, and the half-line defined above otherwise. The angle $\prec_x([x,c],[x,z])$ is defined as the angle between $y_1 -x$ and $z_1 -x$ in the Euclidean space $\Lambda$, where $ y_1 \in [x,c]$ and $z_1 \in [x,z]$ are arbitrary points different from $x$. We endow $\ov{A}$ with the topology generated by the open sets in $A$ and by all of these cones. 

Then $\ov{A}$ is homeomorphic to the ball 
$A_1 = \{x \in A: d(a,x) \leq 1\}$ in $A$. Namely, we can embed $A$ in $A_1$ as 
\begin{eqnarray*}
j: & A & \longrightarrow A_1 \\
~ & x &  \longmapsto \left\{ \begin{array}{ll}
a+ \frac{1 - e^{-d(a,x)}}{d(a,x)} (x-a) & \mbox{ if } x \neq a, \\
a & \mbox{ if } x = a,\end{array} \right.
\end{eqnarray*}
and we map a half-line $c \in A_\infty$ to the point $ x \in c$ with $d(a,x) = 1$,
cf. \cite{schst}, IV.2.
Note that $\ov{A}$ is independent of the choice of $a$.

For every $n \in N$ the affine bijection $\nu(n)$ can be continued to a homeomorphism 
\[\nu(n): \ov{A} \rightarrow \ov{A},\]
since $d$ is $\nu(n)$-invariant. This yields a continuous action of $N$ on $\ov{A}$. 

We will use a description of the Borel-Serre compactification due to Schneider and Stuhler (in \cite{schst}, IV.2) which is formally similar to the definition of $X$. 

Let $\Phi = \Phi(\bt, \bg)$ be the root system corresponding to $\bt$. It consists of finitely many elements in $X^\ast(\bt) \subset \Lambda^\ast$. For all $a \in \Phi$ there  exists a unique closed, connected, unipotent subgroup $\bu_a$ of $\bg$ which is normalized by $\bt$ and has Lie algebra $g_a= \{ X \in g: \Ad(t)X = a(t) X \mbox{ for all }t \in T\}$ (see \cite{bo}, 21.9). We denote the $K$-rational points of $\bu_a$ by $U_a$. 

Now we define for a boundary point $c \in A_\infty$ 
\begin{eqnarray*}
P_c  & = & \mbox{subgroup generated by T and the groups  } U_a \\
~ & ~ & \mbox{ for all } a\in \Phi \mbox{ such that } c \in \{ \ov{x \in A : a(x) \geq 0}\}.
\end{eqnarray*}
Note that if $y \in A_1$ is the point corresponding to $c \in A_\infty$ via the map $j$ defined with $0 \in A$, then $c\in \{ \ov{x \in A : a(x) \geq 0}\}$ iff $a(y) \geq 0$. 
Now we define an equivalence relation on $G \times \overline{A}$ by
\[
(g,x) \sim (h,y)  \mbox { if there exists an } n \in N \mbox{ such that } nx=y \mbox{ and } g\inv h n \in P_x\]
(using the old groups $P_x$ for points $x \in A$). Let $\ov{X}$ be the quotient
\[\ov{X} = G \times \ov{A} / \sim .\]
Then $G$ acts on $\ov{X}$ via left multiplication on the first factor. The compactified appartment $\ov{A}$ can be embedded as $x \mapsto (1,x)$.
Besides, $P_x$ is the stabilizer of $x \in \ov{A}$, and we have a natural $G$-equivariant embedding  $X \hookrightarrow \ov{X}$.

Let $X_\infty = \ov{X} \backslash X$ be the boundary of $X$.
Then $X_\infty$ is the Tits building corresponding to $G$. To be more precise, let $\Delta$ be the simplicial complex whose simplices are the parabolic subgroups of $G$ with  the face relation $P \leq Q $ iff $ Q \subset P$.  Therefore vertices in $\Delta$ correspond to maximal parabolic subgroups $P \subset G$ with $P \neq G$.

Let $| \Delta |$ be the  geometric realization. Then we have a $G$-equivariant bijection 
\[\tau:  | \Delta | \rightarrow X_\infty,\]such that for any $b \in |\Delta|$ the stabilizer of $\tau(b)$ is the parabolic subgroup corresponding to the simplex of $\Delta$ containing $b$ in its interior, see \cite{clt}, 6.1.

There is a natural bijection between parabolic subgroups of $G$ and flags in $V$, associating to a flag in $V$  its stabilizer in $G$. Here maximal proper parabolics correspond to minimal non-trivial flags, hence to non-trivial subspaces $W$ of $V$. 
For any subspace $W$ of $V$ we denote by $y_W$ the vertex in $X_\infty$ corresponding to $W$.

We will now investigate geodesics in $X$, i.e. maps
\[ c : \R \longrightarrow X\]
such that $d(c(t_1), c(t_2)) = |t_1 - t_2|$ for all $t_1, t_2 \in \R$.  A map $c:  \R_{\geq 0} \rightarrow X$ with the same isometry property is called a half-geodesic. 

Note that for any $x_0 \in X$ and for any point $y \in  X _\infty$ there is a unique half-geodesic $\gamma$ in $X$, starting at $x_0$ and converging to $y$ (see \cite{schst}, IV.2). We write $ \gamma = [x_0, y]$. If $x_0$ is also a vertex,  we can describe $\gamma$ as follows: 

\begin{lem}
Fix a vertex $x_0= \{ M\}$ in $X$. Let $W$ be a  non-trivial subspace of $V$ and $y_W$ the corresponding vertex in $X_\infty$. Let $w_1, \ldots, w_n$ be a base of $V$ such that $M = \sum_{i=1}^n R w_i$ and such that $W$ is generated by $w_1,\ldots, w_r$. Then the vertices in $[x_0, y_W]$ are exactly the lattice classes
\[\{ R w_1 + \ldots + R w_r + \pi^k R w_{r+1} + \ldots + \pi^k R w_n\}\]
for all integers $k \geq 0$.
\end{lem}
{\bf Proof: }Note that there exists a basis $w_1, \ldots, w_n$ as in our claim. (By the invariant factor  theorem, we find an $R$-basis $w_1, \ldots,w_n$ of $M$ such that $\alpha_1 w_1, \ldots, \alpha_r w_r$ is an $R$-basis of $M \cap W$ for some $\alpha_i \in K\tim$.)

After applying a suitable  element $g \in PGL(V)$, we can assume that $M= \sum_{i = 1}^n R v_i$ and $ W = \sum_{i=1}^r K v_i$ for our fixed basis $v_1,\ldots, v_n$. Hence $\{M\} = 0 \in \Lambda$. Let $\gamma$ be the half-line
\[ \gamma(t) =c \, t \sum_{i \leq r} \eta_i \quad \mbox{for all } t\geq 0,\]
where $c>0$ is a constant so that $d(0,c \sum_{i \leq r} \eta_i) =1$. We  denote the associated point in $A_\infty$ by $z$.

Now we want to determine $P_z$. We denote by $\chi_i: \bt \rightarrow \Gm$ the character induced by mapping a diagonal matrix to its $i$-th entry. Then $\Delta= \{ a_{i,i+1}= \chi_i - \chi_{i+1}: i = 1, \ldots, n-1\}$ is a base of the root system $\Phi$. By $\Phi^+$ we denote the set of positive roots. If  $a = \sum_i n_i a_{i,i+1}$  is an arbitrary root, then 
\[ z \in  \{\ov{x \in A : a(x) \geq 0}\} \mbox { iff } n_r \geq 0.\]
Hence $P_z$ is generated by $T$ and all groups $U_a$ for all $a = \sum n_i a_{i,i+1}$ with $n_r \geq 0$. The set of roots fullfilling this condition is equal to $\Phi^+ \cup [I]$, where $I = \Delta \backslash \{ a_{r, r+1}\}$ and where $[I]$ denotes the set of roots which are linear combinations of elements in $I$. Hence $P_z$ is the standard parabolic subgroup corresponding to $I$ (see \cite{boti}, 4.2). Therefore it fixes the flag $0 \subset W \subset V$. Hence $z = y_W$, i.e. $\gamma = [0,y_W]$. The vertices on $\gamma$ are the points $k \sum_{i \leq r} \eta_i$ for all integers $k \geq 0$, hence they correspond to the module classes 
\[\{ \pi^{-k}R v_1 + \ldots + \pi^{-k}R v_r + R v_{r+1} + \ldots + R v_n\} = \{ R v_1 + \ldots + R v_r + \pi^{k} R v_{r+1} + \ldots + \pi^{k} R v_n\},\]
as desired.\hfill$\Box$

\begin{lem}
Let $W$ and $W'$ be nontrivial subspaces of $V$. Then the vertices $y_W$ and $y_{W'}$ in $X_\infty$ can be connected by a geodesic in $X$ iff $W \oplus W' = V$. 
\end{lem}
{\bf Proof: }Assume that $y_W $ and $y_{W'}$ can be connected by a geodesic $\gamma$, i.e. $\ga(t) \rightarrow  y_W$ as $t  \rightarrow \infty$, and $\ga(t) \rightarrow y_{W'}$ as $t \rightarrow -\infty$.  By \cite{br}, Theorem 2, p 166, $\gamma$ lies in an appartment. Since our claim is $G$-invariant, we can assume that $\gamma$ lies in our standard appartment $A$. 

After replacing $\ga$ by a parallel geodesic in $A$, we can assume that $\ga$ contains the vertex $x_0 = 0$. Furthermore, after reparametrization we have $\ga(0) = x_0$. So the restriction of $\ga$ to $\R_{\geq 0}$ is equal to $[x_0, y_W]$. 
Let $\Phi^+$ be the set of positive roots corresponding to the base $\Delta= \{a_{1,2}, a_{2,3} \ldots, a_{n-1,n}\}$ of $\Phi$. Let $D$ be the sector 
\[D=\{ x \in A : a(x) \geq 0 \mbox{ for any } a \in \Delta\}\]
in $A$. Since $D$ is a fundamental domain for the operation of the Weyl group, we can furthermore assume that 
$[x_0,y_W] \subset D$. Hence $y_W $ is contained in the boundary of $D$. For every point $z$ in the boundary of $D$ let $I \subset \Delta$ be the set of all $a_{i,i+1}$ such that $[0,z] \subset  \{ x \in A: a_{i,i+1}(x) =0\}$. Then $P_z$ is generated by $T$ and all $U_a$ for $a \in \Phi^+ \cup [I]$, where $[I]$ is the set of roots which are linear combinations of elements in $I$. Hence $P_z$ is the standard parabolic corresponding to $I$. 

Now $P_{y_W}$ is a maximal proper parabolic, hence for $z = y_W$ the set $I$ is just $\Delta \backslash \{a_{r,r+1}\}$ for some $r \leq n-1$. Therefore $W $ is generated by $v_1,\ldots, v_r$. 

By the proof of Lemma 2.2 we know that $[x_0, y_W]$ is the half-geodesic $\ga_1(t) = c_1 t \sum_{i \leq r} \eta_i$ for $t \geq 0$, where $c_1 >0$ is a suitable constant. This half-geodesic can be uniquely continued to a geodesic in $A$, by letting $t$ run over the whole of $\R$.  Since $\gamma$ lies in $A$, we find that $\gamma(t) = c_1 t \sum_{i \leq r} \eta_i$ for all $t \in \R$. 
The proof of Lemma 2.2 also shows that the half-geodesic $\ga_2(t) = c_2 t \sum_{i\leq n-r} \eta_i$ for $t \geq 0$ (and some $c_2 > 0$)  connects $x_0$ with the vertex in $A_\infty$ corresponding to the vector space $K v_1 + \ldots + K v_{n-r}$. Let $p \in N$ be the permutation matrix with
\[ p v_1 = v_{r+1}, \ldots, p v_{n-r} = v_n, p v_{n-r+1} = v_1,\ldots, p v_n = v_r.\]
Then the half-geodesic $p \ga_2$ connects $x_0$ with the point in $A_\infty$ corresponding to the vector space $K v_{r+1} + \ldots+ K v_n$. 
Since $p$ maps $\sum_{i \leq n-r} \eta_i$ to $\sum_{i \leq r} (-\eta_i)$, we have  $p \ga_2(t) = - c_2  t \sum_{i \leq r} \eta_i$, so that $ W' = K v_{r+1} + \ldots + K v_n$, which implies that $W$ and $W'$ are indeed complementary. 

Now assume that $V = W \oplus W'$. Since our claim is $G$-equivariant, we can assume that $W = K v_1 + \ldots +K v_r$ and that $W' = K v_{r+1} + \ldots + K v_n$ for our standard basis $v_1,\ldots, v_n$ and for some $r$ between $1$ and $n-1$. Then

\[\gamma(t) =c_1 t \sum_{i \leq r} \eta_i\]
is a geodesic in $A$ connecting $y_W$ and $y_{W'}$. \hfill$\Box$

Note that this result shows that two vertices in $X_\infty$ can be connected by a geodesic iff the corresponding parabolic subgroups are opposite in the sense of \cite{bo}, 14.20.

We call a geodesic in $X$ combinatorial, if it consists of $1$-simplices and their vertices. The proof of 2.3 shows that any geodesic in $X$ connecting two vertices in $X_\infty$ which contains a vertex in $X$ is already combinatorial. We will now describe the combinatorial geodesics connecting two fixed vertices on the boundary of $X$. 

\begin{prop}
Let $W$ and $W'$ be non-trivial complementary subspaces of $V$, i.e. $W \oplus W' = V$. Let $M$ and $M'$ be lattices of full rank in $W$ respectively $W'$. Then the vertices $\{ M + \pi^k M'\}$ in $X$ for all $k \in \Z$ define a combinatorial geodesic connecting $W$ and $W'$. In fact, this induces a bijection between the set of pairs of lattice classes  $(\{M\},\{M'\})$ with $M \subset W$ and $M' \subset W'$ and the set of combinatorial geodesics connecting $W$ and $W'$ (up to reparametrization). 
\end{prop}
{\bf Proof: }We show first that the vertices $\{M + \pi^k M'\}$ form indeed the vertices of a combinatorial geodesic. After applying a suitable $g \in G$, we can assume that $M = \sum_{i \leq r} R v_i$ and that $M' = \sum_{i \geq r+1} R v_i$ for our standard basis $v_1, \ldots, v_n$ and $r = \dim W$. Then $\gamma(t ) = c_1 t\sum_{i \leq r} \eta_i$ is a combinatorial geodesic in $A$ containing exactly the vertices $\{M + \pi^k M'\}$. (Here $c_1$ is again a constant so that $\ga$ is an isometry). It is clear that up to reparametrization $\ga$ is the unique geodesic in $X$ containing all those vertices.

If $M$ is equivalent to $N$ and $M'$ is equivalent to $N'$, i.e. $M = \pi^a N$ and $M' = \pi^b N'$ for some $a,b \in \Z$, then $M+ \pi^k M' $ is equivalent to $N + \pi^{k + b-a} N'$, hence  the geodesic defined by $M$ and $M'$ coincides with the one defined by $N$ and $N'$ up to reparametrization. 

Assume that $\gamma $ is a combinatorial geodesic connecting $W$ and $W'$, and assume that $\ga(t) { \rightarrow } y_W$ as $t \rightarrow \infty$. As in the proof of 2.3, there is an element $g \in G$ such that $g \gamma$ is contained in our standard appartment $A$. Since $\ga$ is combinatorial, it contains a vertex which we can move to $0 \in A$ by applying some $t \in T$. After reparametrization, we can therefore assume that $ g \ga(0) = 0$. The proof of 2.3 shows furthermore that after composing $g$ with some element of the Weyl group, we can assume that $g \ga|_{\R_{\geq 0}}$ is contained in the sector $D$, which implies that $g \ga(t) = c_1 t \sum_{i \leq r} \eta_i$ for $r = \dim W$. For $M = \sum_{i =1}^r R v_i$ and $M' = \sum_{i = r+1}^n R v_i$ this is the geodesic determined by the vertices $\{M + \pi^k M'\}$ for all $k \in \Z$. Hence $\gamma$ is given by the pair $( g\inv \{ M\},g\inv \{ M'\})$. 

Now suppose that $(\{M\},\{ M'\})$  and $(\{N\},\{ N'\})$ yield the same geodesic $\gamma$. Then, after taking suitable representatives of our module classes, there exists a $k_0 \in \Z$ such that $M+ M' = N + \pi^{k_0}N'$. Put $r= \dim W$, and fix a basis $w_1,\ldots, w_r$ of $W$ such that $M = \sum_{i=1}^r R w_i $ and a basis $w_{r+1}, \ldots, w_n$ of $W' $ with $M' = \sum_{i=r+1}^n R w_i$. Let $A \in GL(r,K)$ and $B \in GL(n-r,K)$ be matrices with $A M = N $ and $B M' = N'$. Then
\begin{eqnarray*}
\left( \begin{array}{cc} A & 0 \\ 0 & \pi^{k_0} B \end{array} \right)
\in GL(n,R),
\end{eqnarray*}
hence $A $ is contained in $GL(r,R)$ and $\pi^{k_0} B$ is contained in $GL(n-r,R)$, which means that $M = N$ and $M' = \pi^{k_0} N'$, i.e. $\{M\} = \{N\}$ and $\{M'\} = \{ N'\}$. \hfill$\Box$

Let $\PP(V) = \mbox{Proj Sym}V^\ast$ be the projective space corresponding to our $n$-dimensional vector space $V$, where $V^\ast$ is the linear dual of $V$. Every non-zero linear subspace $W$ of $V$ defines an integral (i.e. irreducible and reduced) closed subscheme $\PP(W)=\mbox{Proj Sym}W^\ast \hookrightarrow  \PP(V)$ of codimension $n - \dim W$. These cycles given by subspaces of $V$ are called linear. 

Every lattice $M$ (of full rank) in $V$ defines a model $\PP(M) = \mbox{Proj Sym}_R (M^\ast)$ of $\PP(V)$ over $R$, where $M^\ast$ is the $R$-linear dual of $M$. If the lattices $M$ and $N$ differ by multiplication by some $\lambda \in K\tim$ then the corresponding isomorphism $\PP(M) \iso \PP(N)$ induces the identity on the generic fibre.

We call a non-trivial submodule $N$ of $M$ split, if the exact sequence $0 \rightarrow N \rightarrow M \rightarrow M/N \rightarrow 0$ is split, i.e. if $M/N$ is free (or, equivalently, torsion free).
Every split $R$-submodule $N$ of $M$ defines a closed subscheme  $\PP(N)=\mbox{Proj Sym}N^\ast \hookrightarrow \PP(M)$.  It is integral and has codimension $n - \mbox{\rm rk} N$ (see \cite{we}, Lemma 3.1). These cycles in $\PP(M)$ induced by split submodules are also called linear. 

Let $y= y_W$ be a vertex in $X_\infty$ corresponding to the subspace $W$ of $V$, and let $x = \{M\}$ be a vertex in $X$. Then the half-line $[x,y_W]$ connecting $x$ and $y_W$ is combinatorial. We will now show that $[x, y_W]$ governs the reduction of the linear cycle $\ov{\PP(W)}$ induced by $W$ on the model $\PP(M)$.

\begin{prop} Let $M$ be a lattice in $V$, and let 
$x = \{ M\}$ be the corresponding vertex in $X$. For all vertices $y$ in $X_\infty$ let $[x,y]_m$ denote the initial segment of the combinatorial half-geodesic $[x,y]$ consisting of the first $m+1$ vertices. If $y = y_W$, we write $Z_y$ for the linear cycle on $\PP(V)$ defined by $W$. Then we have a bijection
\begin{eqnarray*}
\{ [x,y]_m: y \mbox{ vertex in } X_\infty\} & \longrightarrow & \{ \mbox{ linear cycles in }\PP(M) \otimes_R R / \pi^m\} \\
~ [x,y]_m & \longmapsto & \ov{Z_y} \otimes_R R/\pi^m, 
\end{eqnarray*}
where $\ov{Z_y}$ denotes the closure of $Z_y$ in $\PP(M)$. Hence the initial segments $[x,y_1]_m$ and $[x,y_2]_m$ coincide iff the reductions of $Z_{y_1}$ and $Z_{y_2} $ in $\PP(M) \otimes_R R/\pi^m$ coincide.
\end{prop}

{\bf Proof: }Fix a vertex $y = y_W$ in $X_\infty$. We will first determine the closure $\ov{Z_y}$ in $\PP(M)$. Put $L = W \cap M$. Then $L$ is a free (since torsionfree) $R$-module of rank $r = \dim W$. It is easy to see that the quotient of $L \hookrightarrow M$ is a free $R$-module. Hence $L$ is a split submodule of $M$, so that $\PP(L)$ is an integral closed subscheme of $\PP(M)$.

Obviously, the generic fibre of $\PP(L)$ is equal to $\PP(W) = Z_y$. Hence $\ov{Z_y} =  \PP(L)$.

We define linear cycles on $\PP(M) \otimes_R R/\pi^m$ as cycles $\PP(N) \hookrightarrow \PP(M) \otimes_R R/\pi^m$ for split $R/\pi^m$-submodules $N \hookrightarrow M \otimes_R R/\pi^m$. For such a split submodule $N$ let $N'$ be its preimage in $M$. Note that $N'$ contains $\pi^m M$, so that it has rank $n$. By the invariant factor theorem, we find an $R$-basis $x_1,\ldots, x_n$ of $M$ and non-negative integers $a_1,\ldots, a_n$ such that $\pi^{a_1}x_1, \ldots, \pi^{a_n} x_n$ is a basis of $N'$. Since $N = N' /   \pi^m M$ is a split submodule of $M / \pi^m M$, all $a_i$ must be equal to zero or $m$. We can assume that $a_1 = \ldots = a_r = 0$, and $a_{r+1} = \ldots = a_n = m$. Then $N''= R x_1 + \ldots + R x_r$ is a split submodule of $M$ with reduction $N$, hence $\PP(N)$ is equal to $\PP(N'') \otimes_R R/\pi^m$. Since $\PP(N'')$ is equal to the closure of $Z_{y_W}$ for $W = K x_1 + \ldots + K x_r$, we see that our map is surjective.

We will now show that for split submodules $L_1$ and $L_2$ of $M$ 
\[  \mbox{Proj Sym} L^\ast_1 \otimes_R R/\pi^m = \mbox{Proj Sym} L^\ast_2 \otimes_R R/\pi^m \mbox{ iff } L_1 + \pi^m M = L_2 + \pi^m M.\]

First suppose that $L_1 + \pi^m M = L_2 + \pi^m M$. Then $L_1 \otimes_R R/\pi^m = L_2 \otimes_R R/\pi^m,$ hence by dualizing we find that 
\[ L_1^\ast \otimes_R R/\pi^m = L^\ast_2 \otimes_R R/\pi^m.\]
So $\mbox{Proj Sym}( L_1^\ast \otimes_R R/\pi^m)$ and $\mbox{Proj Sym}( L_2^\ast \otimes_R R/\pi^m)$ coincide as subschemes of $\PP(M) \otimes_R R/\pi^m$, which gives one direction of our claim.

Now assume that $\mbox{Proj Sym}L^\ast_1 \, \otimes_R R/\pi^m = \mbox{Proj Sym}L^\ast_2 \,  \otimes_R R/\pi^m$. Let $\cI_1$ and $\cI_2$ be the corresponding quasi-coherent ideal sheaves on $\PP(M) \otimes_R R/\pi^m$. We denote the quotients of $L_j \hookrightarrow M$ by $Q_j$. Then $Q_j^\ast$ is a free $R$-module, hence $Q_j^\ast \otimes_R R/\pi^m$ is free over $R/ \pi^m$. Let $I_j$ be the ideal in $\mbox{Sym}( M^\ast \otimes_R R/\pi^m)$ generated by a basis of $Q_j^\ast \otimes_R R/\pi^m$. Then $I_j$ coincides with the kernel of the natural map
$\mbox{Sym} ( M^\ast \otimes_R R/\pi^m) \rightarrow \mbox{Sym} ( L_j^\ast \otimes_R R/\pi^m)$. By \cite{ega2}, 2.9.2, we find that $\cI_j = I_j^\sim$.

Therefore $I_1^\sim = I_2^\sim$. Now fix a basis $x_1,\ldots, x_n$ of $M^\ast\otimes_R R/\pi^m$. For every homogeneous ideal $I$ in $S = \mbox{Sym} (M^\ast \otimes_R R/\pi^m)$ let $\mbox{Sat}(I)$ be the ideal in $S$  defined as
\[ \mbox{Sat}(I)= \{ s \in S : \mbox{ for all } i = 1,\ldots, n \mbox{ there exists a } k\geq 0 \mbox{ such that } s x_i^k \in I\}.\]
Then it is easy to check that $I^\sim= J^\sim$ iff $\mbox{Sat} I = \mbox{Sat} J$, compare \cite{ha}, ex. 5.10, p.125. 

Hence we find that $\mbox{Sat} I_1 = \mbox{Sat } I_2$. Since we can choose $x_1,\ldots, x_n$ so that for some $r$ the subset $x_1,\ldots, x_r$ is a basis of $Q_1^\ast \otimes_R R/\pi^m$, it is easy to see that $\mbox{Sat} I_1 = I_1$. Similarly, 
$\mbox{Sat} I_2 = I_2$. Hence $I_1 = I_2$, which implies that $L_1^\ast \otimes_R R/\pi^m$ and $L_2^\ast \otimes_R R/\pi^m$ are isomorphic as quotient modules of $M^\ast \otimes_R R/\pi^m$. Since all three are free over $R/\pi^m $, we find that $L_1 \otimes_R R/\pi^m$ and $L_2 \otimes_R R/\pi^m$ are equal as submodules of $M \otimes_R R/\pi^m$, which implies our claim. 

Now we will show that 
\[ [x,y_{W_1}]_m = [x, y_{W_2}]_m \quad \mbox{ iff} \quad  (W_1 \cap M) + \pi^m M = (W_2 \cap M) + \pi^m M.\]
Let $w_1,\ldots, w_n$ be an $R$-basis of $M$ so that $W_1$ is generated by  $w_1, \ldots, w_r$. Then by Lemma 2.2, the vertices on $[x, y_{W_1}]$ are given by the module classes $\{M_k\}$ for $k \geq 0$, where $M_k=  R w_1 + \ldots + R w_r + \pi^k R w_{r+1} + \ldots+ \pi^k R w_n$. Now  $M_k = (W_1 \cap M ) + \pi^k M$. Hence we find 
\[\begin{array}{cc}
(W_1 \cap M) + \pi^m M = (W_2 \cap M) + \pi^m M & \mbox{ iff } \\
(m+1)\mbox{-th vertex on }[x, y_{W_1}] = (m+1)\mbox{-th vertex on }[x, y_{W_2}]  & \mbox{ iff }\\ ~[x, y_{W_1}]_m = [x, y_{W_2}]_m, & ~
\end{array}\]
as claimed.\hfill$\Box$

This result justifies why we regard $X$ as a kind of graph of $\PP(M)$: its combinatorial geometry keeps track of the reduction of linear cycles.

\section{Non-Archimedean intersections}
Let us fix a lattice $M$ in $\PP(V)$. This defines a smooth, projective model $\Omega= \PP(M)$ of $\PP(V)$ over $R$. By $Z^p(\om)$ we denote the codimension $p$ cycles on $\om$, i.e. the free abelian group on the set of integral (i.e. irreducible and reduced) closed subschemes of codimension $p$. 
If $T \subset \om $ is a closed subset,  we write $CH^p_T(\om)$  for the Chow group of codimension $p$ cycles supported on 
$T$ (see \cite{giso}, 4.1).
For irreducible closed subschemes $Y$ and $Z$ of codimension $p$ respectively $q$ in  $\om$ we can define an intersection class $Y \cdot Z \in 
 CH^{p+q}_{Y \cap Z}(\om)$, 
see \cite{fu}, 20.2 and \cite{giso}, 4.5.1.
We denote by $\deg$ the degree map for 0-cycles in the special fibre of $\om$, i.e. for all $z = \sum n_P P \in Z^d(\om_k)$ we put $\deg z = \sum n_P \; [k(P):k]$, where $k(P)$ is the residue field of $P$.

Let  $Y \in Z^p(\om)$ and $Z \in Z^q(\om)$ be two irreducible closed subschemes such that $p+q=n$ which intersect properly on the generic fibre of $\om$. (Recall that $n$ is the dimension of $V$.) This means that their generic fibres are disjoint, so that $Y \cap Z$ is contained in the special fibre $\om_k$ of $\om$. Hence  we can define a local intersection number
\[<Y,Z> = \deg( Y \cdot Z),\]
where we take the degree of the image of  $Y \cdot Z \in CH^{n}_{Y \cap Z}(\om)$ in $CH^{n-1}(\om_k)$.

We will now fix linear subspaces $A$, $B$, $C$ and $D$ of $V$, such that 
$A$ and $B$ have dimension $p$, and $C$ and $D$ have dimension $q$ for some
$p,q \geq 1$ with $p+q = n$. We will always assume that $q \geq p$. Besides, we assume that the intersections $A \cap C$, $A \cap D$, $B \cap C$ and $B \cap D$ are all zero. 

This implies that the intersection number 
$< \ov{\PP(A)} - \ov{\PP(B)}, \ov{\PP(C)} - \ov{\PP(D)} >$
is defined, where as before $\ov{\PP(A)}$ denotes the closure of the linear cycle $\PP(A)$ in the model $\PP(M)$. 
From Theorem 3.4 in \cite{we} we can deduce that 
\[ < \ov{\PP(A)}, \ov{\PP(C)} > = v ( \det(f_j(a_i))_{i,j= 1,\ldots, p}),\]
where $a_1, \ldots, a_p$ is an $R$-basis of $A \cap M$, and where $f_1, \ldots, f_p$ is an $R$-basis of the free $R$-module $(M / C \cap M)^\ast$. (In other words, $f_1, \ldots, f_p$ are elements in $M^\ast$ generating the ideal corresponding to the linear cycle $\ov{\PP(C)}$.)

Hence we find that 
\[< \ov{\PP(A)} - \ov{\PP(B)}, \ov{\PP(C)} - \ov{\PP(D)} >= v\left( \frac{\det(f_j(a_i)) \det(g_j(b_i))}
{\det(f_j(b_i)) \det(g_j(a_i))}\right) \]
for certain $K$-bases $a_1, \ldots, a_p$ of $A$, $b_1, \ldots, b_p$ of $B$ and $f_1, \ldots, f_p$ of $(V/C)^\ast$, $g_1,\ldots,g_p$ of $(V/D)^\ast$. 
Since the right hand side is invariant under arbitrary base changes of these  vector spaces, the intersection number on the left hand side is independent of the choice of a lattice $M$ in  $V$.
Hence we will also write $< {\PP(A)} - {\PP(B)}, {\PP(C)} - {\PP(D)} >$ for this intersection number.

Now we can prove a formula for such  a local intersection number in terms of the combinatorial geometry of the Bruhat-Tits building $X$. We will call two lattices equivalent, if they define the same lattice class, i.e. if they differ by a factor in $K\tim$.

\begin{thm}
Let $A$, $B$, $C$ and $D$ be as above and assume additionally that $C + D = V$. 
Besides, we assume that there are complementary subspaces $C'$ respectively $D'$ of $C \cap D$ in $C$ respectively $D$, and full rank lattices $L_A$ in $A$ and $L_B$ in $B$ such that the following two conditions hold: 

First of all, the vector space $<A,B>$ generated by $A$ and $B$ is contained in  $C'\oplus D'$.
Secondly, the lattice $p_{C'}(L_A)$ is equivalent to $p_{C'} (L_B)$, and the lattice $p_{D'} (L_A)$ is equivalent to $p_{D'}(L_B)$, where $p_{C'}$ and $p_{D'}$ denote the projections with respect to the decomposition $V = (C \cap D) \oplus C' \oplus D'$.

Choose a lattice $M_0$ in $C \cap D$ and put $M_{C'} = p_{C'}(L_A)$ and $M_{D'} = p_{D'} (L_A)$.
By 2.4 there is a geodesic $\gamma$ corresponding to $M_0 \oplus M_{C'} $ and $M_{D'}$ which connects $C$ and $D'$. We orient $\gamma$ from $C$ to $D'$.  Let $A \ast \gamma$ be the vertex on $\gamma$ closest to the boundary point $y_A \in X_\infty$ in the following sense: $A \ast \gamma$ is the first vertex $x$ on $\gamma$ such that the half-geodesic $[x,y_A]$ intersects $\gamma$ only in $x$. Then
\[<{\PP(A)} - {\PP(B)}, {\PP(C)} - {\PP(D)}  >= p \; \mbox{\rm distor}_\gamma (A\ast \gamma, B\ast \ga),\]
where $\mbox{\rm distor}_\gamma$ means oriented distance along the oriented geodesic $\ga$.
\end{thm}

It would be desirable to rephrase the conditions on $A$, $B$, $C$ and $D$ we impose in 3.1 in terms of the geometry of the boundary $X_\infty$. Before we prove this theorem, let us formulate a corollary in the case $p=1$, where our conditions are rather mild. Note that in this case the intersection pairing we are considering coincides with N\'eron's local height pairing (compare \cite{giso}, 4.3.8).
\begin{cor}
Let $a$ and $b$ be different points in $\PP(V)(K)$ and let $H_C$, $H_D$ be two different hyperplanes in $\PP(V)$ such that the cycles $a-b$ and $H_C -H_D$ in $\PP(V)$ have disjoint supports. Denote by $A$ and  $B$ the lines in $V$ corresponding to $a$ and $b$, and by $C$ and  $D$ the codimension $1$ subspaces in $V$ corresponding to $H_C = \PP(C)$ and $H_D= \PP(D)$.

If $ <A,B> \cap C = <A,B> \cap D$, then $<a-b, H_C-H_D> =0$.

Otherwise, choose lattices $N$ in $C \cap D$, $M_1$ in $<A,B> \cap C$ and $M_2$ in $<A,B> \cap D$. Then $N \oplus M_1 $ is a lattice in $C$, and by 2.4 there is a geodesic $\gamma$ corresponding to $N \oplus M_1 $ and $M_2$ which connects $C$ and $ <A,B> \cap D$. We orient $\gamma$ from $C$ to $<A,B> \cap D$. Let $a\ast \ga$ be the vertex on $\ga$ closest to the boundary point $y_A$. Then  
\[<a-b,H_C-H_D> = \mbox{\rm distor}_\gamma (a\ast \gamma, b\ast \ga).\]
\end{cor}
{\bf Proof of Corollary 3.2: } If the one-dimensional vector spaces $<A,B> \cap C $ and $<A,B> \cap D$ are equal and, say, generated by $w$, we take generators $w_a$ of $A$ and $w_b$ of $B$, and write $w = \lambda w_a + \mu w_b$ with non-zero coefficients $\lambda$ and $\mu$. If $f$ is a homogeneous equation for $C$ and $g$ is a homogeneous equation for $D$, we have $0 = f(w) = \lambda f(w_a) + \mu f(w_b)$ and similarly $0 = \lambda g(w_a) + \mu g(w_b)$. Hence $<a-b, H_C-H_D> = v(f(w_a)/f(w_b))- v(g(w_a)/g(w_b))$ is indeed zero. 

If   $<A,B> \cap C $ and $<A,B> \cap D$ are not equal, we can apply Theorem 3.1.\hfill$\Box$

{\bf Proof of Theorem 3.1: } Note that our conditions imply that $C \cap D$ has dimension $2q-n=q-p$, so that $(C \cap D) \oplus C' \oplus D' = V$.  Since $A$ and $D$ have trivial intersection, $p_{C'}$ induces an isomorphism $p_{C'}|_A: A \rightarrow C'$, so that $M_{C'}$ is indeed a lattice of full rank in $C'$.

Fix an $R$-basis $w_{2p+1}, \ldots, w_n$ of the lattice $M_0$ in $C \cap D$.
Let $a_1, \ldots, a_p$ be an $R$-basis of the lattice $L_A$. If we denote the projection of $a_i$ to ${C'}$ by $w^\sim_{p+i}$, and the projection to $D'$ by $w^\sim_{i}$, we get a basis $w^\sim_1, \dots, w^\sim_p$ of $M_{D'}$ and a basis $w^\sim_{p+1}, \ldots, w^\sim_{2p}$ of $M_{C'}$. Since $A$ is contained in $C' \oplus D'$, we have $a_i = w_i^\sim + w_{p+i}^\sim$. 

Since  $M_{C'} = p_{C'} (L_A)$ is equivalent to $p_{C'}(L_B)$, we can find a constant $\alpha \in K\tim$ such that $p_{C'}(L_B) = \alpha M_{C'}$. Similarly, we find some $\beta \in K\tim$, such that $p_{D'}(L_B) = \beta M_{D'}$. We can write an $R$-basis  $b_1, \ldots, b_p$ of $L_B$ as
\[b_1 = \beta w_1 + \alpha w_{p+1}, \ldots, b_p = \beta w_p + \alpha w_{2p}\]
for some $R$-bases $w_1, \ldots, w_p$ of $M_{D'}$ and $w_{p+1}, \ldots, w_{2p}$ of $M_{C'}$. 

Hence we can calculate the intersection number: 
\[<{\PP(A)} - {\PP(B)}, {\PP(C)} - {\PP(D)}  >=  p\,  v\left(\frac{\alpha} {\beta}\right ).\]

The vertices on $\gamma$ are the lattice classes 
$\{\pi^k \sum_{i=1}^p R w_i + \sum_{i=p+1}^n R w_i\}\quad \mbox{for } k\in \Z$. 
Note that orienting $\gamma$ from $C$ to $D'$ means following these lattice classes in the direction of decreasing $k$.

We will now determine $B \ast \gamma$. Let let $x= \{M\}$ be a vertex on $\ga$, where $M = \pi^k M_{D'} + (M_{C'} + M_0) =  \sum_{i=1}^p \pi^k R  w_i + \sum_{i=p+1}^n R w_i$ for some integer $k$.

Let us first assume that $k > v(\beta/\alpha)$, and put $k_0 = k +v(\alpha) - v(\beta) >0$. Then \[\pi^{k_0} \alpha\inv b_1, \ldots, \pi^{k_0} \alpha\inv b_p, w_{p+1}, \ldots, w_n\] is an  $R$-basis of $M$. Hence by Lemma 2.2 the vertices in $[x,y_B]$ correspond to the module classes $\{\sum_{i=1}^p R \pi^{k_0} \alpha\inv b_i  + \pi^l \sum_{i=p+1}^n  R w_i\}$ for all $l\geq 0$. The vertex next to $x$ on this half-geodesic is 
\[\{ \sum_{i=1}^p R \frac{\pi^{k_0}}{\alpha} b_i + \pi \sum_{i=p+1}^n  R w_i\}= \{\sum_{i=1}^p R \frac{\pi^{k_0 -1}}{ \alpha} b_i  +   \sum_{i=p+1}^n  R w_i\},\]
which lies on $\ga$ since $k_0 -1\geq 0$.
Hence $[x,y_B]$ does not meet  $\ga$ only in $x$. 

Now assume that $k \leq v(\beta/\alpha)$. Then $ \alpha\inv b_1, \ldots, \alpha\inv b_p , \pi^k w_1, \ldots \pi^k w_p, w_{2p+1}, \ldots, w_n$ is an $R$-basis of $M$. Hence by 2.2, the vertices in $[x,y_B]$ correspond to the module classes $\{\sum_{i=1}^p R \alpha\inv b_i + \pi^{k+l} \sum_{i=1}^p R w_i + \pi^l \sum_{i= 2p+1}^n R w_i\}$ for all $l\geq 0$. Therefore the vertex next to  $x$ on this geodesic is $\{ \sum_{i=1}^p R \alpha\inv b_i + \pi^{k+1} \sum_{i=1}^p R w_i + \pi \sum_{i=2 p+1}^n R w_i\}$. Let us assume for the moment that this vertex is contained in $\ga$. Then the module $\sum_{i=1}^p R \alpha\inv b_i + \pi^{k+1} \sum_{i=1}^p R w_i + \pi \sum_{i=2 p+1}^n R w_i$ is equivalent to $\sum_{i=1}^p R w_i + \pi^{-l} \sum_{i=p+1}^n R w_i$ for some $l \in \Z$, which implies that there exists an integer $l_0$ such that the matrix
\[\left( \begin{array}{ccc}
\frac{\beta}{\alpha} \pi^{-l_0} I_{p} & \pi^{k+1-l_0} I_p & 0 \\
\pi^{l -l_0} I_p  & 0 & 0\\
0 & 0 & \pi^{1+l-l_0}I_{q-p}\\
\end{array}
\right)\]
is in $GL(n,R)$. Here $I_p$ denotes the $p \times p$-unit matrix. 

Now we have to distinguish the cases $p\neq q$ and  $p=q$ (where the last block of rows is non-existent). 

Let us first assume that $p \neq q$. Since all entries of this matrix are in $R$ and the determinant is a unit in
$R$, we get $1+l-l_0=0$ and $l-l_0 =0$ which is a contradiction. 

Therefore, in the case $p \neq q$, for $k \leq v(\beta/\alpha)$, the half-geodesic  $[x,y_B]$ meets $\ga$ only in $x$.

If $p = q$, and $k = v(\beta/\alpha)$ we find in a similar way that $l_0 = k+1 =v(\beta/\alpha) + 1$ and $v(\beta/\alpha) - l_0 \geq 0$, which is a contradiction. Hence for this vertex $x$ the half-geodesic $[x,y_B]$ meets $\ga$ only in $x$. 

If however $k < v(\beta/\alpha)$, we find that $\sum_{i=1}^p R \alpha\inv b_i + \pi^{k+1} \sum_{i=1}^p R w_i = \pi^{k+1} \sum_{i=1}^p R w_i +  \sum_{i=p+1}^n R w_i$, which implies that for all these $x$ the half-line $[x,y_B]$ meets $\ga$ not only in $x$. If $p = q$, the vertex $x$ corresponding to $k = v(\beta/\alpha)$ is therefore the only vertex on $\ga$ such that $[x,y_B]$ meets $\ga$ exclusively in $x$. 

In any case we have shown that in our orientation of $\ga$ the vertex $x$ corresponding to $M_B:= (\beta/\alpha) M_{D'}  + M_{C'} + M_0$ is the first vertex on $\ga$ such that $[x,y_B]$ meets $\ga$ only in $x$. Hence $B \ast \ga$ is well defined and given by the module $M_B$. 

In a similar way, we can show that $A  \ast \ga$ is induced by the class of the module $M_A := M_{D'}  + M_{C'} + M_0$. Hence we can calculate
\begin{eqnarray*}
\distor ( A \ast \ga, B \ast \ga)& =& 
\distor (\{M_A\}, \{M_B\}) \\
~ & = &\left\{ \begin{array}{ll}
|v(\beta) - v(\alpha)| ,& \quad \mbox{if }  0 \geq v(\beta) - v(\alpha) \\
-|v(\beta) - v(\alpha)| ,& \quad \mbox{if }  0 < v(\beta) - v(\alpha) \\
\end{array} \right. \\
~ & = & v(\alpha) - v(\beta),
\end{eqnarray*}
which implies our claim. \hfill$\Box$

The proof of Theorem 3.1 also proves the following corollary, which generalizes Manin's formula on $\PP^1$ to higher dimensions (see \cite{ma}, 3.2). 

\begin{cor}
Let $n = 2p$, and let $A$, $B$, $C$ and $D$ be vector spaces in $V$ of dimension $p$, such that $A \oplus C  = A \oplus D = B \oplus C = B \oplus D = C \oplus D = V$. Assume that there are lattices $L_A$ in $A$ and $L_B$ in $B$ such that $p_{C}(L_A)$ is equivalent to $p_{C} (L_B)$, and the lattice $p_{D} (L_A)$ is equivalent to $p_{D}(L_B)$, where $p_{C}$ and $p_{D}$ denote the projections with respect to the decomposition $V =  C \oplus D$.

Put $M_{C} = p_{C}(L_A)$ and $M_{D} = p_{D} (L_A)$.
By 2.4 there is a geodesic $\gamma$ corresponding to $ M_{C} $ and $M_{D}$ which connects $C$ and $D$. We orient $\gamma$ from $C$ to $D$, and we denote by  $A \ast \gamma$  the unique vertex $x$ on $\gamma$ such that the half-geodesic $[x,y_A]$ intersects $\gamma$ only in $x$. Then
\[<{\PP(A)} - {\PP(B)}, {\PP(C)} - {\PP(D)}  >= p \; \mbox{\rm distor}_\gamma (A\ast \gamma, B\ast \ga)\]

\end{cor}

We think of the point $A \ast \ga$ on $\ga$ as the gate to $\gamma$ when entering $\ga$ from the point $y_A$ at infinity. In fact, if $n =2p$, then passing from $y_A$ to a vertex $x$ on $\ga$ means going first to $A \ast \ga$ and then tracking along $\ga$ until $x$ is reached. If $n \neq 2p$, then this holds for all the vertices $x$ on $\ga$ before $A \ast \ga$. 

\section{The symmetric space and its compactification}
Let $Z$ be the symmetric space $G/K$ corresponding to  the real Lie group $G = SL(n,\C)$ and its maximal compact subgroup $K = SU(n,\C)$ for some $n \geq 2$. By ${g} = sl(n,\C)$ and $k = su(n,\C)$ we denote the corresponding Lie algebras. Let $\sigma: G \rightarrow G$ be the involution $A \mapsto (^t \ov{A})\inv$, and put $p = \{X \in g : d\sigma X  =-X\}$. Note that $d\sigma X =- ^t \ov{X}$ and $k = \{ X \in g: d\sigma X = X\}$. We denote by $\Ad_G : G \rightarrow GL(g)$ the adjoint representation. Then $p$ is invariant under $\Ad_G K$ and $g = p +k$. 

Let $\tau:G \rightarrow G/K = Z$ be the projection map, mapping $1$ to $u \in Z$. The homomorphism $d\tau: g \rightarrow T_uZ$ induces an isomorphism $d \tau : p \simeq T_u Z$  with $d\tau (\Ad(k) X) = d \lambda(k) d\tau(X)$ for all $k \in K$ and $X \in p$, where $\lambda(g)$ denotes the left action of $g\in G$ on $G/K$. (See \cite{he}, IV,\S 3.)

Let $B : g \times g \rightarrow \R$ defined by 
\[B(X,Y) = \Tr(\ad X \ad Y) = 4n \mbox{Re}\, \Tr (XY)\]
  be the Killing form on $g$, see \cite{he}, III, 6.1 and III, 8. $B$ is positive definite on $p \times p$ and induces a scalar product $<,>$ on $T_u Z$ via the isomorphism $d\tau: p \simeq T_u Z$. Shifting this product around with the $G$-action, we get a $G$-invariant metric on $Z$. We write $\dist(x,y)$ for the corresponding distance between two points in $Z$. For any $X \in p$, the geodesic in $u$ with tangent vector $d\tau X$ is given by 
\[\ga(t) = (\exp tX) u,\]
where $\exp: g\rightarrow G$ is the exponential map, induced by the matrix exponential function (see \cite{he}, IV,3).
The geodesic connecting two points in $Z$ can be described as follows:
\begin{lem}
Let $z_1$ and $z_2$ be two points in $Z$. Then there exists an element $f \in G$ such that $ z_1 = f u$ and $ z_2 = f d u$, where $d$ is a diagonal matrix of determinant one with positive real entries $d_1,\ldots, d_n$. Put $a_i = \log d_i \in \R$, and let $X$ be the diagonal matrix with entries $a_1,\ldots, a_n$. The geodesic connecting $z_1$ and $z_2$ is
\[\ga(t)= f \exp (t X) u \quad \mbox{for } t \in [0,1], \]
and we have $ \dist(z_1,z_2) = 2 \sqrt{n} \sqrt{\sum a_i^2}$.
\end{lem}
{\bf Proof: } A straightforward calculation.\hfill$\Box$

We will now describe the differential geometric compactification of $Z$ using half-geodesics (see e.g. \cite{bgs} \S 3 or \cite{jo}, 6.5).
A ray emanating at $z \in Z$ is a unit speed (half-)geodesic $\ga: \R_{\geq 0} \rightarrow Z$ with $\ga(0) = z$. Two rays $\ga_1$ and $\ga_2$ are called asymptotic if $\dist(\ga_1(t), \ga_2(t))$ is bounded in $t \in \R_{\geq 0}$. The set of equivalence classes of rays with respect to this relation is denoted by $Z(\infty)$, and we put $\ov{Z} = Z \cup Z(\infty)$.

\begin{lem}
For any $z \in Z$ and any $c \in Z(\infty)$ there is a unique ray $\ga$ starting in $z$ whose class is $c$. We refer to $\gamma$ as the geodesic connecting $z$ and $c$.\end{lem}
{\bf Proof: } See \cite{jo}, Lemma 6.5.2, p. 255.\hfill$\Box$

This fact implies that for each $z \in Z$ we can identify $Z(\infty)$ with the unit sphere $S_z Z = \{ X \in T_zZ: ||X||=1\}$ in $T_zZ$ by associating to an element $X \in S_z Z$ the equivalence class of the ray $\ga$ with $\gamma(0) = z$ and $\stackrel{\cdot}{\gamma}(0) = X$.

For $z \in Z$ and $c_1,c_2 \in Z \cup Z(\infty)$ such that $ z\neq c_1$ and $z \neq c_2$ we define $\prec_z(c_1,c_2)$ as the angle between $\stackrel{\cdot}{\ga_1}(0)$ and $\stackrel{\cdot}{\ga_2}(0)$ in $T_z Z$, where $\ga_1$ and $\ga_2$ are the geodesics from $z$ to $c_1$ respectively $c_2$. For $z \in Z$, $c \in Z(\infty)$ and $\epsilon >0$ let $C_z(c,\epsilon)$ be the cone
\[ C_z(c,\epsilon) = \{ y \in \ov{Z}: y \neq z \mbox{ and }\prec_z(c,y) < \epsilon\} .\]
The cone topology on $\ov{Z}$ is the topology generated by the open sets in $Z$ and these cones, see \cite{bgs}, 3.2,  p.22.

The bijection $S_z Z \rightarrow Z(\infty)$ is then a homeomorphism. Since every $g \in G$ acts by isometries on $Z$, it acts in a natural way on $Z(\infty)$, and the corresponding action on $\ov{Z}$ is continuous (\cite{bgs}, 3.2, p.22).

A flat $E$ in $Z$ is a complete totally geodesic Euclidean submanifold of maximal dimension (which is by definition the rank of $Z$). Let $E$ be a flat in $Z$ with $u \in E$. Then $T_u E$ is a maximal abelian subalgebra in $T_u Z \simeq p$ and $E =\tau(\exp a)$. In fact, the flats in $Z$ containing $u$ correspond bijectively to the maximal abelian subspaces of $p$, see \cite{jo}, 6.4.2, p.248.

Let $a$ be a maximal abelian subspace of $p$. It is easy to see that the bilinear form on $g$ defined as 
\[<X,Y>_g = \left\{ \begin{array}{ll}
B(X,Y) & \mbox{if } X,Y \in p \\
- B(X,Y) & \mbox{if } X,Y \in k \\
0 & \mbox{ if } X \in p, Y \in k \mbox{ or } X \in k, Y \in p 
\end{array} \right. \]
is positive definite and that $ \{ \ad H: H \in a\}$ is a commuting family of self-adjoint endomorphisms with respect to $<,>_g$  (compare \cite{jo}, \S 6.4). Hence $g$ can be decomposed as an orthogonal sum of common eigenspaces of the $\ad H$:
\[g = g_0 \oplus \bigoplus_{\lambda \in \Lambda} g_\lambda,\]
where
\begin{eqnarray*}
g_\lambda = \{ X \in g: \ad(H) (X) = \lambda(H) X \mbox{ for all } H \in a\} \mbox{, and }\\
g_0 = \{ X \in g: \ad(H) (X) = 0 \mbox{ for all } H \in a\} \mbox{, and where}\\
\Lambda \mbox{  is the set of all  }\lambda \neq 0 \mbox{  in  } \mbox{Hom}_\R ( a, \R) \mbox{ such that  }g_\lambda \neq 0.
\end{eqnarray*}
Note that $\Lambda$ is a root system  (see \cite{kn}, Corollary 6.53).

Let $E$ be a flat containing $u$.
A geodesic $\ga: \R \rightarrow E$ with $\ga(0) = u$ is called singular if it is also contained in other flats besides $E$; if not, it is called regular. Tangent vectors to regular (singular) geodesics are also called regular (respectively singular). A vector $H \in a$ is singular  iff there is an element $Y \in g \backslash g_0$ such that $[H,Y] = 0$, hence iff there exists a $\lambda \in \Lambda$ with $\lambda(H) =0$, see \cite{jo}, 6.4.7, p.251.
We denote by $a_{\rm sing}$ the subset of singular elements in $a$, i.e. $a_{\rm sing} = \{H \in a: \mbox{ there exists some } \lambda \in \Lambda \mbox{ with } \lambda(H)=0\}$, and by $a_{\rm reg}$ the complement $a_{\rm reg} = \{ H \in a: \mbox{ for all } \lambda \in \Lambda: \lambda (H) \neq 0\}$. The singular hyperplanes $\cH_\lambda = \{ H \in a: \lambda (H) = 0\}$ divide $a_{\rm reg}$ into finitely many components, the Weyl chambers.

More generally, a face of $a$ with respect to the $\cH_\lambda$ is defined as an equivalence class of points in $a$ with respect to the following equivalence relation: $x \sim y$ iff for each $\cH_\lambda$, $x$ and $y$ are either both contained in $\cH_\lambda$ or lie on the same side of $\cH_\lambda$. The faces not contained in any hyperplane $\cH_\lambda$ are called chambers (see \cite{bou}, p. 60). The set of chambers in $a$ corresponds to the set of bases of the root system $\Lambda$ in the following way: If $B$ is a basis of $\Lambda$, then the corresponding chamber $C = C(B)$ can be described as
\[C = \{ X \in a: \lambda (X) > 0 \mbox{ for all  } \lambda \in B\}.\]
Besides, the faces contained in $\ov{C}$ correspond bijectively to the subsets of $B$, if we associate to $I \subset B $ the set
\[ C_I = \{  X \in a : \lambda(X) = 0 \mbox{ for all } \lambda \in I \mbox{ and } \lambda(X) >0 \mbox{ for all } \lambda \in B\backslash I\}.\]
(See \cite{bou}, V.1 and Th\'eor\`eme 2 in VI, 1.5.) Note that every face in $a$ is contained in the closure of a chamber (\cite{bou}, Proposition 6, p. 61.)
We will only consider the faces of dimension bigger than zero, i.e. we will always assume that $ I \neq B$. 

We can now transfer the faces in $p$ to the boundary: 
For every maximal abelian subspace $a \subset p$ and every face $F \subset a$ we denote by $F(\infty)$ the so-called face at $\infty$
\[F(\infty) = \{ [\ga]: \ga(t) = \exp (t H)u   \mbox{ for some } H \in F \mbox{ of norm }1\} \subset Z(\infty),\]
where $[\ga]$ is the equivalence class of the ray $\ga$. We will now  investigate these faces at infinity. 

Note that $G$ can be identified with the set of $\R$-rational points of a semisimple linear algebraic group $\bg$ over $\R$ (which can be defined as a suitable subgroup of  the algebraic group $SL(2n, \R)$). For every closed algebraic subgroup of $\bg$ the group of $\R$-rational points is a Lie subgroup of $G$, since it is $\R$-closed  (see \cite{wa}, Theorem 3.42).

\begin{lem} All maximal abelian subspaces of $p$ are Lie algebras of maximal $\R$-split tori in $\bg$. If $a$ and $a'$ are maximal abelian subspaces of $p$, then there exists an element $k \in K$ with $\Ad(k) a = a'$.  
\end{lem}
{\bf Proof: }Any abelian subspace $a$ of $p$ consists of pairwise commuting hermitian matrices over $\C$, hence there exists an orthonormal basis of $\C^n$ of common eigenvectors of the elements in $a$.
Therefore for some $k \in K$ the maximal abelian subspace $k a k\inv$ is contained in (and hence equal to) the  abelian subspace $d \subset p$ of real diagonal matrices with trace $0$. Since $d = \Lie T$, where $T = \bt(\R)$ for the maximal $\R$-split torus $\bt$ in $\bg$ of real diagonal matrices with determinant $1$, our claim follows.\hfill$\Box$

Let now $\bt$ be a maximal $\R$-split torus in $\bg$ such that  $a =\Lie T$ (for T = $\bt(\R)$) is a maximal abelian subspace of $p$, and let $g = g^T \oplus \bigoplus_{\alpha \in \Phi(T,G)} g_\alpha$ be the root decomposition corresponding to $\bt$ (see \cite {bo}, 8.17 and 21.1). Here $g^T = \{ X \in g: \Ad(t) X = X \mbox { for all } t \in T\}$, and $g_\alpha = \{X \in g: \Ad(t) X = \alpha(t) X \mbox { for all } t \in T\}$, and $\Phi(\bt,\bg)$ is the set of roots, i.e. the set of non-trivial characters of $\bt$ such that $g_\alpha \neq \{0\}$. 

This is the same decomposition as the one we defined previously for the maximal abelian subspace $a$ of $p$. Namely, by choosing a basis of  $g^T$ and all $g_\alpha$, we find a basis of $g$ such that $\Ad(t) \in \mbox{GL}(g)$ is given by a diagonal matrix for all $t \in T$. Passing to the Lie algebras, we find that 
$g^T  = g_0$ and $g_\alpha = g_\lambda $ for $\lambda = d \alpha \in a^\ast$. In particular, we have an additive bijection
\[  \Phi(\bt,\bg) \longrightarrow \Lambda,\]
which we will use from now on to identify $\Phi(\bt,\bg)$ and $\Lambda$. 

Let $B$ be a base of the root system $\Phi(\bt,\bg)$, and let $I \subset B$ be a proper subset. By $[I]$ we denote the set of roots which are linear combinations of elements in $I$, and we set $\psi(I) = \Phi^+ \backslash [I]$, where $\Phi^+$ denotes the set of positive roots with respect to $B$. Let $\bu_{\psi(I)}$ be the unique closed connected unipotent subgroup of $\bg$, normalized by $Z(\bt)$, with Lie algebra $\bigoplus_{\alpha \in \psi(I)} g_\alpha$ (\cite{bo}, 21.9), and let $\bt_I$ be the connected component of the intersection $\cap_{\alpha \in I} \ker \alpha$. Then the standard parabolic subgroup $\bp_I$ corresponding to $I$ is defined as the semidirect product $\bp_I = Z(\bt_I) \bu_{\psi(I)}$. Now $\Lie(Z(\bt_I)) = g^T \oplus \bigoplus _{\alpha \in [I]} g_\alpha$. Hence $\Lie \bp_I = g^T \oplus \bigoplus_{\alpha \in [I] \cup \psi(I)} g_\alpha$. Note that we exclude the trivial parabolic $\bg$ here since we do not allow $I = B$. 

Every proper parabolic subgroup of $\bg$ is conjugate to a uniquely determined standard parabolic by an element in $G = \bg(\R)$  (\cite{bo}, 21.12). 

Let us put $\bp = \bp_\emptyset$. Then $\bp$ is a minimal parabolic in $\bg$. Let $\bn$ be the normalizer of $\bt$ in $\bg$, and put $N = \bn(\R)$, $P = \bp(\R)$. Then $(G,P,N)$ gives rise to a Tits system by \cite{bo}, 21.15, i.e. it fulfills the conditions in $\cite{bo}, 14.15$. 

Hence $(P,N)$ is a $BN$-pair for $G$ in the terminology of \cite{br}, V.2. A subgroup $Q$ of the group $G$ is called parabolic if $Q$ contains a conjugate of $P$. As this terminology suggests, we have a bijection
\[ \bq \mapsto Q = \bq(\R)\]
between the set of parabolic subgroups of $\bg$ and the set of parabolic subgroups of $G$ (see \cite{bo}, 21.16). The Tits building $Y$ corresponding to the $BN$-pair $(P,N)$ is defined as the partially ordered set (poset) of proper parabolic subgroups of $G$ with by the relation $Q_1 \leq Q_2 $ if $Q_2 \subset Q_1$. This poset is in fact the poset of simplices of a simplicial complex (\cite{br}, V.3).

There is a natural $G$-equivariant bijection 
between proper parabolic subgroups of $G$ and non-trivial flags in $\C^n$, associating to a flag its stabilizer in $G$, so that we can identify $Y$ with the poset of flags in $\C^n$.

We can now describe the stabilizers of points in $Z(\infty)$ as follows:

\begin{prop}
Let $a$ be a maximal abelian subspace of $p$, let $C$ be a chamber in $a$ and $B$ the associated base of the root system $\Lambda$. Let $X$ be a vector in the face $C_I \subset \ov{C}$ for some $I \subset B$ satisfying $||X||=1$. Then we denote by $\ga(t) = \exp (tX) u$ the ray in $u$ defined by $X$, and by $z $ the corresponding point in $Z(\infty)$. 
Let $G_z$ be the stabilizer of $z$. Then $G_z$ is the standard parabolic subgroup $P_I= \bp_I(\R)$ corresponding to $I$.
\end{prop}
{\bf Proof: } As in  \cite{jo}, Theorem 6.2.3, one can show  that $\Lie G_z = g_0 \oplus \bigoplus_{\lambda (X) \geq 0} g_\lambda$. Since $X$ is in $C_I$, we have $\lambda(X) =0$ for all $\lambda \in I$ and $\lambda(X) >0$ for all $\lambda \in B\backslash I$. We denote by $\Phi^+$ respectively $\Phi^-$ the positive respectively negative roots with respect to $B$. Then $\lambda(X) \geq 0$ for all $\lambda \in \Phi^+$. Besides,  $\lambda (X) \geq 0$ for some $\lambda \in \Phi^-$ iff $\lambda \in [I]$. Hence $\Lie G_z = g^T \oplus \bigoplus_{\alpha \in \Phi^+ \cup [I]} g_\alpha = \Lie P_I$. Note that the minimal parabolic  $P_\emptyset$ stabilizes $z$, so that it is contained in $G_z$. Hence $G_z$ is a standard parabolic, and we must have  $G_z = P_I$.\hfill$\Box$

The following corollary shows that the poset of faces at infinity is isomorphic to the Tits building $Y$ of $G = SL(n,\C)$.
\begin{cor}[cf. \cite{bgs}, p.248f]
The association $F \longmapsto G_z$ for any $z \in F(\infty)$ defines an order preserving bijection between the set of faces (of dimension $>0$)   in all the maximal abelian subspaces of $p$ (ordered by the relation $F' < F $ iff $F' \subset \ov{F} $) and the set of proper parabolic subgroups of  $G$ ordered by the relation $P' < P$ iff $P \subset P'$.
\end{cor}
{\bf Proof: }Let us first assume that two faces $F$ and $F'$ are mapped to the same parabolic subgroup. Hence there are elements $X \in F$ and $X' \in F'$  of norm 1 such that the points $z$ and $z' $ in $Z(\infty)$ corresponding to the geodesics $(\exp tX) u$ and  $(\exp t X' )u$, respectively, satisfy $G_z = G_{z'}$. Choose maximal abelian subspaces $a$ and $a'$ of $p$ containing $F$ respectively $F'$ with corresponding root systems $\Lambda$ and $\Lambda'$. There is a base $B$ of $\Lambda$ and a subset $I \subset B$ such that $F = C_I$  for the chamber $C$ in $a$ given by $B$, and, similarly, a base $B'$ of $\Lambda'$ and a subset $I' \subset B'$ such that $F' =  C'_{I'}$, where $C'$ is the chamber in $a'$ induced by $B'$. 

By Lemma 4.3 there exists an element $k \in K $ such that $\Ad(k)a' = a$. It is easy to see that the map 
\begin{eqnarray*}
a^{'\ast} & \stackrel{\varphi}{\longrightarrow}& a^\ast \\
\lambda' & \longmapsto & \lambda' \circ \Ad(k\inv)
\end{eqnarray*}
induces a bijection beween $\Lambda'$ and $\Lambda$ so that
\[ \Ad(k) g_{\lambda'} = g_{\varphi(\lambda')}.\]
Hence the homomorphism $\Ad(k): a' \rightarrow a$ maps hyperplanes to hyperplanes and thus faces to faces.

So  $\Ad(k) C'_{I'}$ is a face in $a$, and therefore contained in the closure of a chamber in $a$. Let $N_K(a) = \{ k \in K : \Ad(k) a = a \}$ and $Z_K(a) = \{ k \in K : \Ad(k) H = H \mbox{ for all } H \in a\}$ be the normalizer respectively the centralizer of $a$ in $K$. Then $W = N_K(a) / Z_K(a)$ is the Weyl group of $\Lambda$ (see \cite{kn}, 6.57), hence it acts transitively on the set of chambers. 

So we may choose $k$ so that $\Ad(k) C'_{I'}$ is a face in the closure of $C$, i.e. $\Ad(k) C'_{I'} = C_J$ for some $J \subset B$. Now we have for all $t \in \R$
\[ \exp (t \Ad(k) X') u = k \exp (t X') u,\]
hence $kz'$ is the class of the geodesic $(\exp t \Ad(k) X') u$. Applying Proposition 4.4 we find that 
\[ P_J = G_{k z'} = k G_{z'} k \inv = k G_z k\inv = k P_I k\inv,\]
which implies that $I = J$ and that $k \in P_I= G_z = G_{z'}$. Hence we have $kz' = z'$. Since $\exp (t X' )u$ and $\exp (t \Ad(k) X')  u $ are unit speed geodesics connecting $u$ with $z' = k z'$, we have $X' = \Ad(k) X'$ by 4.2. Therefore $X'$ must be in $C_J = C_I$, so that $X$ and $X'$ both  lie in the face $C_I$. 

Let now $Y$ be an arbitrary element of $C'_{I'}$. Then $Y_0 = Y / ||Y||$ is also contained in $C'_{I'}$ and induces a point $z_0 $ in $Z(\infty)$. Since $z_0$ lies in the same face at infinity as $z'$, we have $G_{z_0} = G_{z'}$. Hence the same reasoning as above implies that $Y_0$ and hence $Y$ lies in $C_I$. Altogether we find that $C'_{I'} \subset C_I$, i.e. that $F' \subset F$.  Reversing the roles of $F$ and $F'$, we also have the opposite inclusion, so that  $F = F'$, which proves injectivity. 

This implies  in particular that two faces in $p$ are either disjoint or equal, and also that two faces at infinity are either disjoint or equal. 

Let us now show surjectivity. Fix a chamber $C$ in a maximal abelian subspace $a$ of $p$. A proper parabolic subgroup $P \subset G$ is conjugate to a standard parabolic subgroup associated to $C$, hence, using Proposition 4.4, it is the stabilizer of some $z \in Z(\infty)$. Let $X \in p$ be the unit vector such that $z$ is the class of $(\exp tX) u$, and let $a'$ be any maximal abelian subspace containing $X$. Then $X$ lies in some face $F$ in $a'$, which is mapped to $G_z = P$. 

Now assume that $F$ and $F'$ are two faces in $p$ satisfying $F' \subset \ov{F}$. Let $a$ be a maximal abelian subspace of $p$ containing $F$. Then there exists a chamber $C$ in $a$ such that $F = C_I$ for some subset $I$ of the base corresponding to $C$. Since $F'$ is contained in $\ov{C}$, it meets a face $C_J$ for some $J \subset B$. Since we have already seen that faces are disjoint or equal, we find that $F' = C_J$, so that $I \subset J$. Hence $P_I \subset P_J$, which by 4.4 implies $G_{z} \subset G_{z'}$ for any two points $z \in F(\infty)$ and $z' \in F'(\infty)$. 

On the other hand, assume that $G_{z} \subset G_{z'}$ for two points $z$ and $z'$ in $Z(\infty)$ such that $z$ is the class of $(\exp tX) u$ and $z'$ is the class of $(\exp t X' )u$ for unit vectors $X$, $X'$ in $p$. Take a maximal abelian subspace $a$ and a chamber $C$ such that $X \in C_I$ for some subset $I$ of the base corresponding to $C$. Then $G_{z} = P_I$ by 4.4. Therefore $G_{z'}$ is a standard parabolic, hence we find a set $J \supset I$ with $G_{z'} = P_J$. By injectivity, $X'$ lies in $C_J$, so that $I \subset J$ implies $C_J \subset \ov{C_I}$. Hence  
our map is order preserving in both directions.\hfill$\Box$

The minimal faces of positive dimension in some maximal abelian subspace $a \subset p$ are the faces 
\[ F = \{ H \in a : \lambda_0(H) >0 \mbox{ and } \lambda(H) = 0 \mbox{ for all }H \in B\backslash \{ \lambda_0\}\},\]
where $B$ is a base of $\Lambda$ and $\lambda_0$ is an element of $B$. The corresponding face at  infinity  consists of one point. According to Corollary 4.5, the minimal faces correspond to the maximal proper parabolic subgroups of $G$, which in turn correspond to minimal flags in $G$, i.e. to non-trivial subspaces $ W \subset \C^n$. Hence we find  that the set of all non-trivial subspaces of $\C^n$ can be regarded as  a subset of $Z(\infty)$. 
We denote the point in $Z(\infty)$ corresponding to a subspace $W \subset \C^n$ by $z_W$. We endow $\C^n$ with the canonical scalar product with respect to the standard basis $e_1,\ldots, e_n$. 

\begin{lem} Let $W$ be an $r$-dimensional subspace of $\C^n$ with $ 0 < r <n$. 
Choose an orthonormal basis $w_1,\ldots, w_r$ of $W$ and complete it to an orthonormal basis $w_1,\ldots, w_n$ of $\C^n$ so that the matrix $g$ mapping $e_i$ to $w_i$ for all $ i = 1,\ldots, n$ is contained in $K$. Then $z_W$ is the class of the following ray emanating at $u$: 
\begin{eqnarray*}
\{ g \,\exp t {\ \left(  \begin{array}{cccccc}
\rho & & & & & \\
& \ddots& & & & \\
& & \rho & & & \\
& & &  \sigma & & \\
& & & & \ddots & \\
& && & &  \sigma
\end{array} \right)u}: t \geq 0 \},
\end{eqnarray*}

where $\rho = \sqrt{\frac{n-r}{4rn^2}}$ and  $\sigma = - \frac{r}{n-r} \rho$, and where $\rho$ appears $r$ times.
\end{lem}

{\bf Proof: }Let $\bt$ be the maximal split torus consisting of all real diagonal matrices  of determinant $1$ in $\bg$ with respect to $e_1,\ldots, e_n$, and put $T = \bt(\R)$. 
The corresponding root system in $a = \Lie T$ is 
\[ \Lambda= \{ \lambda_{ij}: i \neq j\}, \]
where $\lambda_{ij} = \lambda_i - \lambda_j$, and $\lambda_i \in a^\ast$ maps a diagonal matrix to its $i$-th entry. The subset $B = \{\lambda_{12}, \ldots , \lambda_{n-1 n}\}$ is a base of $\Lambda$. Let $C$ be the corresponding chamber. 

The vector space $W^\sim = \C e_1 \oplus \ldots \oplus \C e_r$ corresponds to the standard parabolic 
\begin{eqnarray*} P = \left\{ \left( \begin{array}{cc}
{\ast}& \ast \\
0 & \ast
\end{array} \right) 
\begin{array}{cc}
\} r & ~ \\
~& ~
\end{array}\right\} \subset G
\end{eqnarray*}
given by $B \backslash \{ \lambda_{r r+1} \}$. Now put $\rho = \left( \frac{n-r}{4 r n^2}\right)^{1/2}$ and $\sigma = -\frac{r}{n-r} \rho$. Then the diagonal matrix $X$ with entries $({\rho, \ldots, \rho}_r, \sigma, \ldots, \sigma)$ has norm $1$ and is contained in $C_{B \backslash \{ \lambda_{r r+1}\} }$. Hence $z_{W^\sim}$ is given by the ray $\{ \exp (tX)u: t \geq 0 \}$.  Applying $g$, our claim follows.~~~$\quad$~~~~\hfill$\Box$

From now on, we will write
\begin{eqnarray*}
\diag(d_1,\ldots, d_n) = { \left( \begin{array}{ccc}
d_1& & 0 \\
& \ddots& \\
0& & d_n \\
\end{array} \right).}
\end{eqnarray*}
The preceding Lemma says that 
\begin{eqnarray*}
\gamma(t) = g \exp t \diag(\rho, \ldots, \rho ,\sigma, \ldots, \sigma) \cdot u
\end{eqnarray*}
is the ray connecting $u$ and $z_W$. We write $[u, z_W]$ for this ray. 
\begin{cor}
Let $x \in Z$ be an arbitrary point, and let $W$ be an $r$-dimensional  subspace of $\C^n$ with $0 < r <n$.  Put $W^\sim = \C e_1 \oplus \ldots \oplus \C e_r$. Then there exists an element $g \in G $ such that $g u = x$ and $g W^\sim = W$. For any such $g$ let 
\begin{eqnarray*}
\gamma(t) = g \exp t \mbox{\rm diag}(\underbrace{\rho, \ldots, \rho}_r, \sigma, \ldots, \sigma)\cdot u \mbox{ for } t \geq 0,
\end{eqnarray*}
where $\rho = \sqrt{ \frac{n-r}{4 r n^2}}$ and  $\sigma= -\frac{r}{n-r} \rho$. Then $\gamma = [ x, z_W]$. 
\end{cor}
{\bf Proof: }The last assertion is an immediate consequence of Lemma 4.6. Note that an element $g \in G$ as in our claim exists. Namely, let $f \in G$ be any element satisfying $u = fx$ and choose an orthonormal basis $w_1, \ldots, w_r$ of $f(W)$. Then we can complete it to an orthonormal basis $w_1,\ldots, w_n$ of $\C^n$ such that the element $k$ mapping $e_i$ to $w_i$ is contained in $K $. Hence $g = f\inv k$ maps $u$ to $x$ and $W^\sim $ to $W$.~~~\hfill$\Box$

We will now investigate full geodesics in $Z$ connecting two $0$-simplices on the boundary $Z(\infty)$. The following result is the Archimedean analogue of 2.3. 
\begin{lem}
Let $W$ and $W'$ be two non-trivial subspaces of $\C^n$. Then there exists a geodesic joining $z_W$ and $z_{W'}$ iff $W\oplus W' = \C^n$.
\end{lem}
{\bf Proof: } Assume first that $W\oplus W' = \C^n$ and that $\dim W = r$. Applying a suitable $g \in G$, we can assume that $W= \C e_1 + \ldots + \C e_r$  and $W' = \C e_{r+1} + \ldots \C e_n$. Put again $\rho = \sqrt{\frac{n-r}{4 rn^2}}$ and $\sigma= - \frac{r}{n-r} \rho$, and let  
\[\gamma(t) = \exp t \;\diag (\underbrace{\rho,\ldots, \rho}_r, \sigma,\ldots, \sigma) u\]
for all $t \in \R$. For $t \geq 0$, this is equal to $[u, z_W]$ by 4.6. Let $N$ be the permutation matrix mapping $(e_1,\ldots, e_n)$ to $(e_{r+1}, \ldots, e_n, e_1,\ldots, e_r)$. Then we have for all $t \geq 0$
\begin{eqnarray*}
\gamma(-t) & =  &N \exp t  (\underbrace{-\sigma,\ldots,-\sigma}_{n-r},-\rho,\ldots, -\rho) N\inv u\\
& = & N \exp t  (-\sigma,\ldots,-\sigma,-\rho,\ldots, -\rho) u.
\end{eqnarray*}
Since $-\sigma = \sqrt{\frac{r}{4 (n-r)n^2}}$ and $-\rho = - \frac{n-r}{r} (-\sigma)$, we can apply  Lemma 4.6 and find that $ \gamma(-t)$ is the ray $ [ u, z_{W'}]$. Hence $\ga$ is a geodesic in $Z$ connecting $W$ and $W'$.

Now suppose that there exists a geodesic $\ga$ joining $W$ and $W'$, and let $x = \ga(0)$. For $t \geq 0$, the half-geodesic $\ga(t)$ connects $x$ and one of the vector spaces, say $W$. Hence by 4.7 (after a reparametrization of $\gamma$ so that it has unit speed), 
\[\ga(t) = g \exp t \;\diag (\rho,\ldots, \rho, \sigma,\ldots, \sigma) u \quad \mbox{ for } t \geq0,\]
where $g$ maps $\C e_1 + \ldots, \C e_r$ to $W$ and $u$ to $x$. Then this equality holds for all $t \in \R$. We have already seen in the other half of our proof that 
\[\exp(-t) \diag(\rho, \ldots, \rho,\sigma,\ldots, \sigma) u\]
connects $u$ with $\C e_{r+1} + \ldots + \C e_n$, so that $\gamma(-t)$ for $t \geq 0$ connects $x$ with $g(\C e_{r+1} + \ldots + \C e_n)$. Hence $W' = g(\C e_{r+1} + \ldots + \C e_n)$, which implies that $W \oplus W' = V$. \hfill$\Box$

Note that this result - as its non-Archimedean counterpart 2.3 - shows that two minimal faces in $Z(\infty)$ can be connected by a geodesic iff the corresponding parabolic subgroups are opposite. 

Note that the map
\[ Z \ni z = g K \longmapsto g ~^t \ov{g}\]
provides a bijection between points in $Z$ and positive definite hermitian matrices in $SL(n,\C)$  of determinant one, or, what amounts to the same, equivalence classes $\{h\}$ of hermitian metrics, i.e. positive definite hermitian forms, on $\C^n$ with respect to the relation $h \sim h'$, if $h$ is a positive real multiple of $h'$. 

We can now prove an Archimedean analogue of Proposition 2.4. 

\begin{prop}
Let $W$ and $W'$ be complementary subspaces of $\C^n$, i.e. $W \oplus W' = \C^n$, and put  $r = \dim W$. Let $\{h\}$ and $\{h'\}$ be equivalence classes of hermitian metrics on $W$ respectively $W'$, and let $g$ be an element in $SL(n,\C)$ such that $ge_1,\ldots, g e_r$ is an orthonormal  basis of $\alpha h$ and $g e_{r+1}, \ldots, g e_n$ is an orthonormal basis of $\alpha' h'$ for some representatives $\alpha h$ of $\{h\}$ and $\alpha' h'$ of $\{h'\}$. Then 
\[g \;\exp t \;\mbox{\rm diag}(\underbrace{\rho,\ldots, \rho}_r,\sigma,\ldots,\sigma)u \]
is a geodesic connecting $W$ and $W'$, where $\rho= \sqrt{\frac{n-r}{4rn^2}}$ and $\sigma= -\frac{r}{n-r} \rho$.  In fact, this association defines a bijection between the set of pairs
$(\{h\},\{h'\})$ of metric classes on $W$ and $W'$ and the set of geodesics (up to reparametrization) connecting $W$ and $W'$.
\end{prop}
{\bf Proof: }Obviously, given $\{h\}$ and $\{h'\}$, we can always find an element $g$ as in our claim. It is clear that 
$g  \exp t \diag  (\rho,\ldots,\rho,\sigma,\ldots,\sigma) u$ is a geodesic connecting $W$ and $W'$. Let  $f \in SL(n,\C)$ be another element such that $f e_1,\ldots , f e_r$ is an orthonormal basis of $\beta h$ and $f e_{r+1}, \ldots , f e_n$ is an orthonormal basis of $\beta' h'$ for some positive real numbers $\beta$ and $\beta'$. Then there are positive real numbers $\delta$ and $\delta'$ with $\delta^r \delta^{' n-r}=1$ such that $\diag(\delta\inv,\ldots, \delta\inv,\delta^{' -1},\ldots, \delta^{' -1}) g\inv f$ is in $K$. Besides, $g\inv f$ is of the form
\begin{eqnarray*}
g\inv f = \left( \begin{array}{cc}
{\ast}& 0 \\
0 & \ast
\end{array} \right)
\begin{array}{cc}
\} r & ~ \\
~& ~
\end{array}.
\end{eqnarray*}
Hence 
\begin{eqnarray*}
\lefteqn{ f \,\exp t\, \diag\,(\rho,\ldots, \rho,\sigma,\ldots,\sigma)u }\\
& = & g \,  \exp t \,\diag (\rho,\ldots, \rho,\sigma,\ldots,\sigma)(g\inv f) \,u 
\\
&= & g\, \exp t\, \diag (\rho,\ldots, \rho,\sigma,\ldots,\sigma) \diag(\delta,\ldots, \delta,\delta' ,\ldots, \delta') u \\
& = & g \, \exp(t + t_0) \diag(\rho,\ldots, \rho, \sigma,\ldots,\sigma)u,
\end{eqnarray*}
where $t_0$ satisfies $\rho t_0 = \log \delta$ (hence $\sigma t_0 = -\frac{r}{n-r} \rho t_0 = \log \delta'$). This shows that up to reparametrization our geodesic is independent of the choice of $g$. 

Now let $\ga$ be a geodesic connecting $W$ and $W'$. We have seen in the proof of 4.8 that up to reparametrization $\ga(t) = g \exp t \diag (\rho,\ldots,\rho,\sigma,\ldots,\sigma) u$, where $g$ maps $\C e_1 \oplus \ldots \oplus \C e_r$ to $W$ and $\C e_{r+1} \oplus \ldots \oplus \C e_n$ to $W'$. Let $w_i= g(e_i)$ for all $i = 1,\ldots, n$, and let $h$ respectively $h'$ be the metrics on $W$ respectively $W'$ with orthonormal bases $w_1,\ldots, w_r$ respectively $w_{r+1}, \ldots, w_n$. Then $\ga$ is induced by the pair $(\{h\},\{h'\})$. 

Now suppose that $(\{h_1\},\{ h'_1\})$ and $(\{h_2\},\{ h'_2\})$ are pairs of metric classes  leading to the same geodesic. Let $g_1$ and $g_2$ elements in $SL(n,\C)$ such that for $i = 1$ or $2$ $g_i e_1,\ldots, g_i e_r$ is an orthonormal basis of $\alpha_i h_i$ and $g_i e_{r+1}, \ldots, g_i e_n$ is an orthonormal basis of $\alpha_i' h_i'$, where $\alpha_i, \alpha_i'$ are positive real numbers. Then there exists some $t_0 \in \R$ such that 
\[ g_1 u = g_2 \exp t_0 \, \diag(\rho,\ldots, \rho,\sigma,\ldots,\sigma) u.\]
Put $d =  \exp t_0\, \diag(\rho,\ldots,\rho,\sigma,\ldots,\sigma)$. Then $g_1\inv g_2 d$ is contained in $K$. Let us denote by $\lambda_0$ the canonical scalar product on $\C^n$, and let $\lambda_i$ be the metric on $\C^n$ with orthonormal basis $g_i e_1, \ldots, g_i e_n$ for $i = 1$ or $2$. Then $\lambda_i$ is the orthogonal sum of $\alpha_i h_i$ and $\alpha'_i h'_i$. By definition, we have for all $v,w \in \C^n$
\[ \lambda_0(v,w) = \lambda_1(g_1 v, g_1 w) \quad \mbox{and}\quad \lambda_0(v,w) = \lambda_2 (g_2 v, g_2 w).\]
Since $g_1\inv g_2 d$ is in $K$, we find that
\[ \lambda_2 (g_2 v, g_2 w)= \lambda_0(v,w) = \lambda_0 (g_1\inv g_2 d v, g_1\inv g_2 d w).\]
If $v$ and $w$ are in $\C e_1 + \ldots + \C e_r$, then $d v = \delta v$ and $d w = \delta w$ for the positive real number $\delta = \exp( t_0 \rho)$. Hence
\begin{eqnarray*}
\lefteqn{\alpha_2 h_2( g_2 v, g_2 w)  }\\
& = & \lambda_2(g_2 v, g_2 w) = |\delta|^2 \lambda_0(g_1\inv g_2 v , g_1\inv g_2 w) \\
& =&   |\delta|^2 \lambda_1(g_2 v, g_2 w) = |\delta|^2  \alpha_1 h_1 ( g_2 v ,g_2 w),
\end{eqnarray*}
which implies that $h_1$ is equivalent to $h_2$ on $W$. Similarly, looking at vectors in $\C e_{r+1} + \ldots + \C e_n$, we find that $h_1'$ is equivalent to $h_2'$.\hfill$\Box$

\section{Archimedean intersections}

In this section we  prove an Archimedean analogue of Theorem 3.1.
We will first define the local Archimedean intersection number of linear cycles in $\PP^{n-1}_{\C}$. Fix a hermitian metric $h$ on $V = \C^n$. Then we can also define a metric $h^\ast$ on the dual vector space $V^\ast$. Let $W$ be a linear subspace of $V$ of codimension $p$ and let $z_1,\ldots, z_n$ be an orthonormal basis of $V^\ast$ such that $W$ is the intersection $W = \cap_{i=1}^p \mbox{\rm ker} (z_i)$. (Then the linear cycle $\PP(W) \subset \PP(V) = \PP^{n-1}_{\C}$ is given by the homogeneous ideal generated by $z_1,\ldots, z_p$.) On $\PP(V) \backslash \PP(W)$ we define
\[\tau = \log(|z_1|^2 + \ldots + |z_n|^2) \quad \mbox{and} \quad \sigma = \log(|z_1|^2 + \ldots + |z_p|^2).\]
Beside, we define $(1,1)$-forms $\alpha = d d^c \tau$ and $\beta = d d^c \sigma$ on $\PP(V) \backslash \PP(W)$, where $d d^c = \frac{i}{2 \pi} \partial \ov{\partial}$. Put
\[\Lambda_W = (\tau - \sigma) (\sum_{\nu = 0}^{p-1}\alpha^\nu \beta^{p-1-\nu}),\]
which is the Levine form for the linear cycle $\PP(W)$ (see \cite{giso}, 1.4 and \cite{giso2}). Then $\Lambda_W$ induces  a Green current for $\PP(W)$ with associated form $\alpha^p$ (\cite{giso2}, 5.1).

From now on we fix linear subspaces $A$, $B$, $C$ and $D$ of $V=\C^n$, such that 
$A$ and $B$ have dimension $p$, and $C$ and $D$ have dimension $q$ for some
$p,q \geq 1$ with $p+q = n$. We will always assume that $q \geq p$. Besides, we assume that the intersections $A \cap C$, $A \cap D$, $B \cap C$ and $B \cap D$ are all zero. In this case the cycles $\PP(A) - \PP(B)$ and $\PP(C) - \PP(D)$ meet properly on $\PP(V)$. We define their Archimedean intersection number as
\[<\PP(A) - \PP(B), \PP(C) - \PP(D)> = \int_{\PP(A)} (\Lambda_C - \Lambda_D) - \int_{\PP(B)} (\Lambda_C - \Lambda_D),\]
(compare \cite{giso}, 4.3.8, iii).
Using the commutativity of the $\ast$-product for Green currents (\cite{giso}, 2.2.9), we find that this is independent of the choice of a hermitian metric $h$ on $V$.

Now we can prove a formula for such  a local intersection number in terms of the geometry of the symmetric space $Z$. Taking into account the correspondence between lattices on the non-Archimedean side and hermitian metrics on the Archimedean side (cf. \cite{de}), the following result is the Archimedean counterpart of our non-Archimedean Theorem 3.1.

\begin{thm}
Let $A$, $B$, $C$ and $D$ be as above and assume additionally that $C + D = V$. 
Besides, we assume that there are complements $C'$ respectively $D'$ of $C \cap D$ in $C$ respectively $D$, and hermitian metrics $h_A$ on $A$ and $h_B$ on $B$ such that the following two conditions hold: 

First of all, the vector space  $<A,B>$ generated by $A$ and $B$ is contained in  $C'\oplus D'$.
Secondly, the metric $p_{C'\ast}(h_A)$ is equivalent to $p_{C' \ast} (h_B)$, and the metric $p_{D'\ast} (h_A)$ is equivalent to $p_{D'\ast}(h_B)$, where $p_{C'}$ and $p_{D'}$ denote the projections with respect to the decomposition $V = (C \cap D) \oplus C' \oplus D'$.

Choose a metric $h_0$ in $C \cap D$ and put $h_{C'} = p_{C'\ast}(h_A)$ and $h_{D'} = p_{D' \ast} (h_A)$.
By 4.9 there is a geodesic $\gamma$ corresponding to the orthogonal sum $h_0 \oplus h_{C'} $ on $C$ and the metric $h_{D'}$ on $D'$ which connects $C$ and $D'$. We orient $\gamma$ from $C$ to $D'$.  Let $A \ast \gamma$ be the unique point $z$ on $\gamma$ such that the ray $[z, z_A]$ meets $\gamma$ at a right angle. Then 
\[<{\PP(A)} - {\PP(B)}, {\PP(C)} - {\PP(D)}  >= \frac{\sqrt{p}}{\sqrt{q}} \; \mbox{\rm distor}_\gamma (A\ast \gamma, B\ast \ga),\]
where $\mbox{\rm distor}_\gamma$ means oriented distance along the oriented geodesic $\ga$.
\end{thm}

Before we prove this theorem, let us formulate a corollary in the case $p=1$, where our conditions are rather mild. Note that in this case the intersection pairing we are considering coincides with N\'eron's local height pairing (compare \cite{giso}, 4.3.8).

The following result generalizes Manin's formula for $\PP^1$ to higher dimensions (see \cite{ma}, Theorem 2.3).
\begin{cor}
Let $a$ and  $b$ be different points in $\PP^{n-1}_{\C}$, and let $H_C$ and $H_D$ be different hyperplanes in $\PP^{n-1}_{\C}$ such that the cycles $a-b$ and $H_C-H_D$ have disjoint supports. Denote by $A$ and $B$ the lines in $\C^n$ corresponding to $a$ and $b$, and by $C$ and $D$ the codimension $1$ subspaces in $\C^n$ corresponding to $H_C$ and $H_D$.

If $ <A,B> \cap C = <A,B> \cap D$, then $<a-b, H_C-H_D> =0$.

If this is not the case, choose hermitian metrics $h_0$ on $C \cap D$, $h_1$ on $<A,B> \cap C$  and $h_2$ on $<A,B> \cap D$. Then $h_0 \oplus h_1$ is a metric on $C$, and by  4.9 there exists a geodesic $\gamma$ connecting $C$ and $<A,B> \cap D$ associated to the pair $(\{h_0 \oplus h_1\},\{h_2\})$. We orient $\gamma$ from  $C$ to $<A,B> \cap D$. 
Then 
\[<a-b, H_C -H_D> =  \frac{1 }{\sqrt{n-1}} \;\mbox{\rm distor}_\gamma (a\ast \gamma, b\ast \ga).\]

\end{cor}

{\bf Proof of Corrolary 5.2: } If the one-dimensional vector spaces $<A,B> \cap C$ and $<A,B> \cap D$ are equal, a similar reasoning as in the proof of Corollary 3.2 shows that $<a-b, H_C-H_D>$ is indeed zero. If they are not equal, we can apply Theorem 5.1.\hfill$\Box$

{\bf Proof of Theorem 5.1: }Our conditions imply that $C \cap D$ has dimension $2q-n = q-p$, so that indeed $V = (C \cap D ) \oplus C' \oplus D'$. Since $A$ and $D$ have trivial intersection, the projection $p_{C'}: A \rightarrow C'$ is a linear isomorphism, so that we can define $p_{C' \ast} h_A$ as $h_A \circ p_{C'}\inv$.  Hence the orthogonal sum $h_0 \oplus h_{C'}$ is indeed a hermitian metric on $C$. Let now $w_{2p+1},\ldots, w_n$ be an orthonormal basis of $h_0$ in $C \cap D$. 
Let $a_1, \ldots, a_p$ be an orthonormal basis of $h_A$ in $A$. If we denote the projection of $a_i$ to $C'$ by $w_{p+i}^\sim$, and the projection to $D'$ by $w_i^\sim$, we get an orthonormal basis $w_1^\sim, \ldots, w_p^\sim$ of $h_{D'}$, and an orthonormal basis $w_{p+1}^\sim, \ldots, w_{2p}^\sim$ of $h_{C'}$. Since $A$ is contained in $C' \oplus D'$, we have $a_i = w_i^\sim + w_{p+i}^\sim$. 

Since $h_{C'} =  p_{C'\ast}(h_A)$ is equivalent to $p_{C' \ast} (h_B)$, we can find a constant $\alpha \in \R_{>0}$ so that $\alpha^2 p_{C' \ast} (h_B) = h_{C'}$. Similarly, we find some $\beta \in \R_{>0}$ such that $\beta^2 p_{D' \ast}(h_B) = h_{D'}$. Hence there is an orthonormal basis $b_1, \ldots, b_p$ of $h_B$ in $B$ such that 
\[b_1 = \beta w_1 + \alpha w_{p+1}, \ldots, b_p = \beta w_p + \alpha w_{2p}\]
for some orthonormal bases $w_1, \ldots w_p$ of $h_{D'}$ and $w_{p+1}, \ldots, w_{2p}$ of $h_{C'}$.

We define a metric $h$ on $V$ as the orthogonal sum of $h_0, h_{C'}$ and $h_{D'}$. If $\Lambda_C$ and $\Lambda_D$ are the Levine currents with respect to $h$, we can calculate \[\int_{\PP(B)} \Lambda_C - \Lambda_D = 2 p \, \log (\alpha/\beta) \int_{\PP(B)} (d d^c \log (|u_1|^2 + \ldots+ |u_p|^2))^{p-1},\] where $u_1,\ldots, u_p$ are projective coordinates on $\PP(B)$. Since $\frac{i}{2}  \partial \ov{\partial} \log  (|u_1|^2 + \ldots+ |u_p|^2)$ is the $(1,1)$-form with respect to the Fubini-Study-metric on $\PP(B)$, we have \[\int_{\PP(B)} ( \frac{i}{2}  \partial \ov{\partial} \log  (|u_1|^2 + \ldots+ |u_p|^2))^{p-1} = (p-1)! \, \mbox{vol} (\PP(B))\] by Wirtinger's theorem (see \cite{grha}, p. 31). Therefore \[\frac{1}{\pi^{p-1}}  \int_{\PP(B)} \left( \frac{i}{2}  \partial \ov{\partial} \log  (|u_1|^2 + \ldots+ |u_p|^2)\right)^{p-1}  = \frac{(p-1)!}{\pi^{p-1}} \, \mbox{vol} (\PP(B)) = 1,\] 
since the volume of $\PP(B)$ with respect to the Fubini-Study-metric is $\frac{\pi^{p-1}}{(p-1)!}$ (see e.g. \cite{bgm}, p. 18).
A similar calculation gives $\int_{\PP(A)} \Lambda_C - \Lambda_D = 0$, so that  our intersection number is 
\[<{\PP(A)} - {\PP(B)}, {\PP(C)} - {\PP(D)}  >= 2p \log \frac{\beta}{\alpha} .\]

For some complex number $\delta$ the element $g$ mapping 
$e_1,\ldots, e_n$ to $\delta\inv  w_1,\ldots,\delta\inv w_n$ is in $SL(n,\C)$. Then by 4.9, putting $\rho= \sqrt{q/(4 p n^2)}$ and $\sigma = - \frac{p}{q} \rho$,
\[\gamma(t) = g \, \exp t \,\diag (\underbrace{\rho,\ldots, \rho}_p,\underbrace{\sigma,\ldots,\sigma}_q) u\]
is the geodesic connecting $C$ and $D'$ corresponding to $\{h_0 \oplus h_{C'}\}$ and $\{h_{D'}\}$. Orienting $\gamma$ from $C$ to $D'$ means following the direction of  increasing $t$. 

We will now determine $B \ast \ga$. Let $z = g \, \exp t \, \diag(\rho,\ldots,\rho,\sigma,\ldots,\sigma) u$ be a point on $\gamma$. We put $g^\sim = g \, \exp t \, \diag (\rho,\ldots, \rho,\sigma, \ldots, \sigma)$, so that $g^\sim u = z$. We want to describe $[z,z_B]$. Let  $c$ be the positive real number
\[ c = (\beta^2 \exp({-2t \rho}) + \alpha^2 \exp({-2t \sigma} )) ^{\frac{1}{2}},\]
and let $k \in SL(n,\C)$ be the matrix
\begin{eqnarray*}
k = \left(
\begin{array}{ccc}
\frac{\beta}{c} \exp(-t\rho) I_p & - \frac{\alpha}{c}\exp(-t\sigma)I_p& 0 \\
\frac{\alpha}{c}\exp(-t\sigma)I_p& \frac{\beta}{c} \exp(-t\rho) I_p& 0 \\
0 & 0 & I_{q-p}
\end{array}
\right),
\end{eqnarray*}
where $I_p$ denotes the $(p \times p)$-unit matrix. Obviously, $k$ is an element in $K$ such that  $g^\sim k = g \exp t \diag(\rho, \ldots, \rho, \sigma,\ldots,\sigma) k$ maps the vector space generated by $e_1,\ldots, e_p$ to $B$. Using 4.7 we find that $[z,z_B]$ is given by 
\[ g^\sim  \, k \, \exp (s\, \diag(\rho,\ldots, \rho,\sigma,\ldots,\sigma))\,  u \quad \mbox{ for } s \geq 0.\]
On the other hand, $[z,z_{D'}]$ is the ray
\[  g^\sim \, \exp (s \, \diag(\rho,\ldots, \rho,\sigma,\ldots,\sigma)) \,u \quad \mbox{ for } s\geq 0.\]
Now the angle between $[z,z_A]$ and $[z, z_{D'}]$ in $z =    g^\sim \, u$ is equal to the angle between 
\begin{eqnarray*}\gamma_1(s) &=&
k \,\exp (s \; \diag( \rho, \ldots, \rho,\sigma, \ldots,
\sigma) )u \quad \mbox{ and} \\
\gamma_2(s)& =& \exp( s \;  \diag( \rho, \ldots, \rho, \sigma, \ldots, \sigma) )u
\end{eqnarray*}
in $u$. 

Recall that $\tau: G \rightarrow Z$ is the projection map, and that $\lambda(g)$ denotes the left action of $g \in G$ on $Z$. Then we have 
\begin{eqnarray*}
\stackrel{\cdot}{\gamma}_2(0) & = & d \tau \, \diag ( \rho,\ldots, \rho, \sigma,\ldots,\sigma) \quad \mbox{and}\\
\stackrel{\cdot}{\gamma}_1(0) & = & d \lambda(k) d\tau \, \diag(\rho,\ldots, \rho, \sigma,\ldots,\sigma) \\
~&=& d\tau ( \Ad(k) \, \diag(\rho,\ldots, \rho,\sigma,\ldots,\sigma)),\end{eqnarray*}
see section 4. Recall that $<,>$ is the scalar product on $T_u Z$ induced by the Killing form on $p$. We have
\begin{eqnarray*}
\lefteqn{<\stackrel{\cdot}{\gamma}_1(0) , \stackrel{\cdot}{\gamma}_2(0) >
= }\\
& & 4n\mbox{Re} \Tr (k \diag(\rho,\ldots, \rho,\sigma,\ldots,\sigma) k\inv \diag(\rho,\ldots, \rho,\sigma,\ldots,\sigma)).
\end{eqnarray*}
Calculating this matrix, we find that 
\begin{eqnarray*}
\lefteqn{ \Tr (k \diag(\rho,\ldots, \rho,\sigma,\ldots,\sigma) k\inv \diag(\rho,\ldots,\sigma,\ldots,\sigma))}\\
& &= p \left[\frac{\rho^2 \beta^2}{c^2} \exp(-2t\rho)+ 2 \frac{ \alpha^2\rho \sigma}{c^2} \exp(-2t \sigma)+  \frac{\sigma^2 \beta^2}{c^2} \exp(-2t\rho)\right]\; +  (q-p) \sigma^2\\
& &=\frac{\rho^2}{c^2 q^2}\left( p(p^2 + q^2) \beta^2 \exp(-2t\rho) -2p^2q \alpha^2 \exp(-2t\sigma) +(q-p) p^2 c^2 \right)\\
& &= \frac{\rho^2}{c^2 q^2}\left( pqn \beta^2 \exp(-2t\rho)  - p^2 n \alpha^2 \exp(-2t\sigma)\right).
\end{eqnarray*}
Therefore
\begin{eqnarray*}
<\stackrel{\cdot}{\gamma}_1(0) , \stackrel{\cdot}{\gamma}_2(0) > 
= 0,& \mbox{iff} &   \exp(-2t\rho+ 2t\sigma) =  \frac{ p \, \alpha^2}{q \, \beta^2} ,\\
~ & \mbox{hence iff}&  \frac{n}{q} 2t\rho  = \log\left(\frac{q \,\beta^2}{p \, \alpha^2}\right).
\end{eqnarray*}
Thus $[z,z_B]$ meets $[z, z_{D'}]$ (and hence $\gamma$) at a right angle, iff $t = \frac{q}{2n\rho}$ \mbox{$\log(q \beta^2/ p \alpha^2)$}. So $B \ast \gamma $ is well-defined and equal to the point

\[ B \ast \gamma = g \, \diag \left(\left(\frac{q \,\beta^2}{p \, \alpha^2}\right)^{\frac{q}{2n}}, \ldots, \left(\frac{q \,\beta^2}{p \, \alpha^2}\right)^{\frac{q}{2n}}, \left(\frac{q \,\beta^2}{p \, \alpha^2}\right)^{\frac{-p}{2n}}, \ldots,  \left(\frac{q \,\beta^2}{p \, \alpha^2}\right)^{\frac{-p}{2n}}\right) u.\]
An analogous calculation gives 
\[A \ast \gamma  =g  \, \diag \left(\left(\frac{q }{p }\right)^{\frac{q}{2n}}, \ldots, \left(\frac{q}{p}\right)^{\frac{q}{2n}}, \left(\frac{q}{p }\right)^{\frac{-p}{2n}}, \ldots,  \left(\frac{q}{p}\right)^{\frac{-p}{2n}}\right) u.\]
Now we can  calculate
\begin{eqnarray*}
\lefteqn{\dist(A \ast \gamma, B\ast \gamma)}\\
& & = \dist \left( u,\diag\left( \left(\frac{\beta}{\alpha}\right)^{\frac{q}{n}}, \ldots,  \left(\frac{\beta}{\alpha}\right)^{\frac{q}{n}},\left(\frac{\beta}{\alpha}\right)^{\frac{-p}{n}}, \ldots
,\left(\frac{\beta}{\alpha}\right)^{\frac{-p}{n}} \right)u\right)\\
& & = 2\sqrt{pq} \left| \log \frac{\beta}{\alpha}\right|,
\end{eqnarray*}
by Lemma 4.1. 

Now we have $1 \leq \beta/\alpha$ iff $A \ast \gamma$ appears before $B \ast \gamma$ in our orientation of $\gamma$. Hence 
\begin{eqnarray*}
\distor_\gamma ( A\ast \ga, B \ast \ga)& =& 
\left\{ \begin{array}{ll}
2 \sqrt{pq} \left| \log \frac{\beta}{\alpha} \right| & \quad \mbox{if }\alpha \leq \beta\\
-2 \sqrt{pq} \left| \log \frac{\beta}{\alpha} \right| & \quad \mbox{if }\alpha >\beta
\end{array} \right.\\
& = &2 \sqrt {pq}\log \frac{\beta}{\alpha},
\end{eqnarray*}
which implies that
\[<\PP(A) - \PP(B), \PP(C)- \PP(D)>=   2p \log \frac{\beta}{\alpha}  = \frac{\sqrt{p}}{ \sqrt{q}} \; \distor_\ga(A \ast \gamma, B \ast \gamma),\]
whence our claim.\hfill$\Box$

\end{document}